\documentclass{svjour3}                     
\smartqed  
\usepackage{subfigure}
\usepackage{amsmath}
\usepackage{amssymb}
\usepackage{algorithm}
\usepackage{psfrag,pict2e,picture}
\usepackage{epic}
\usepackage{soul}
\usepackage{tikz}

\newcommand{\ba}{\begin{array}}
\newcommand{\ea}{\end{array}}
\newcommand{\bea}{\begin{eqnarray}}
\newcommand{\eea}{\end{eqnarray}}
\newcommand{\be}{\begin{equation}}
\newcommand{\ee}{\end{equation}}
\newcommand{\bd}{\begin{displaymath}}
\newcommand{\ed}{\end{displaymath}}
\newcommand{\bi}{\begin{itemize}}
\newcommand{\ei}{\end{itemize}}
\newcommand{\bn}{\begin{enumerate}}
\newcommand{\en}{\end{enumerate}}
\newcommand{\pa}{\partial}
\newcommand{\f}{\frac}

\newcommand{\hf}{\frac12}

\newtheorem{rem}{Remark}
\newcommand{\bb}[1]{{{#1}}}
\newcommand{\rr}[1]{{{#1}}}

\begin{document}


\title{Hermite -- Discontinuous Galerkin Overset Grid Methods \\
for the Scalar Wave Equation}
\titlerunning{H--DG Overset grid methods}        
\author{Oleksii Beznosov        \and
Daniel Appel\"{o}
}
\institute{Oleksii Beznosov \at
  Department of Mathematics and Statistics, \\
  University of New Mexico, \\
  1 University of New Mexico, MSC01 1115, Albuquerque, NM 87131 \\
              \email{obeznosov@unm.edu}
           \and
Daniel Appel\"{o} \at Department of Applied Mathematics, \\
University of Colorado,\\
University of Colorado 526 UCB, Boulder, CO 80309 \\
\email{daniel.appelo@colorado.edu}
}

\date{Received: date / Accepted: date}



\maketitle

\begin{abstract}
We present high order accurate numerical methods for the wave equation that combines efficient Hermite methods with geometrically flexible discontinuous Galerkin methods by using overset grids. Near boundaries we use thin boundary fitted curvilinear grids and in the volume we use Cartesian grids so that the computational complexity of the solvers approach a structured Cartesian Hermite method. Unlike many other overset methods we do not need to add artificial dissipation but we find that the built in dissipation of the Hermite and discontinuous Galerkin methods is sufficient to maintain stability. By numerical experiments we  demonstrate the stability, accuracy, efficiency and applicability of the methods to forward and inverse problems.
\end{abstract}

\section{Introduction}
Accurate and efficient simulation of waves is important in many areas in science and engineering due to the ability of waves to carry information over large distances. This ability stems from the fact that waves do not change shape in free space. On the other hand when the background medium is changing this induces a change in the wave forms that propagate through the medium and the waves can be used for probing the interior material properties of objects.

In order to preserve the properties of waves from the continuous setting it is preferable to use high order accurate discretizations that are able to control dispersive errors. The development of high order methods for wave propagation problems has been an active area of research for a long time and there are by now many attractive methods. Examples include (but are not limited to) finite difference methods, \cite{SBP,Virta2014,Wang2016,Petersson2018,Hagstrom2012}, embedded boundary finite differences, \cite{appelo2012fourth,FCAD1,FCAD2,Li2004295,Wandzura2004763}, element based methods like discontinuous Galerkin (DG) methods, \cite{Wilcox:2010uq,Upwind2,ChouShuXing2014,ChungEngquist06,ChungEngquist09,GSSwave}, hybridized discontinuous Galerkin (HDG) methods, \cite{Nguyen2011,Stanglmeier2016}, cut-cell finite elements \cite{STICKO2016364,sticko2016higher} and Galerkin-difference methods \cite{BANKS2016310}.

An advantage of summation-by-parts finite differences and Galerkin type methods is that stability is guaranteed, however this guarantee also comes with some drawbacks. For diagonal norm summation-by-parts finite differences the order of accuracy is reduced to roughly half of that in the interior near boundaries. Further the need for multi-block grids also restricts the geometrical flexibility.

As DG and HDG methods are naturally formulated on unstructured grids they have good geometric flexibility.  However, Galerkin based polynomial methods often have the drawback that they require small timesteps (the difference Galerkin and  cut-cell finite element methods are less affected by this) when combined with explicit timestepping methods, but on the other hand they preserve high order accuracy all the way up to the boundary and it is easy to implement boundary conditions independent of the order of the method.

The pioneering work by Henshaw and co-authors, see for example \cite{chess1990}, describe techniques for generating overset grids as well as how they can be used to solve elliptic and first order time-dependent partial differential equations (PDE) by second order accurate finite differences.
In an overset grid method the geometry is  discretized by narrow body-fitted curvilinear grids while the volume is discretized on one or more Cartesian grids. The generation of such body-fitted grids is local and typically produces grids of very high quality, \cite{OGEN}. The grids overlap (we say that they are overset) so that the solution on an interior (often referred to as non-physical or ghost) boundary can be transferred from the interior of another grid. In \cite{chess1990} and in most other overset grid methods the transfer of solutions between grids is done by  interpolation.
 Since the bulk of the domain can be discretized on a Cartesian grid the efficiency asymptotically approaches that of a Cartesian solver but still retains the geometrical flexibility of an unstructured grid method.

\rr{We note that the same type of efficiency can be expected for embedded boundary and cut-cell finite elements. A difference is that overset grid methods typically have smoother errors near physical boundaries and this may be important if quantities that include derivatives of the solution, such as traction or strain, are needed. }

Here we are concerned with the approximation of the scalar wave equation on overset grids. To our knowledge, high order overset grid methods for wave equations in second order form have been restricted to finite difference discretizations. For example, in \cite{henshaw:1730} high order centered finite difference approximations to Maxwell's equations (written as a system of second order wave equations) was introduced. More recently, in \cite{ANGEL2018534}, the upwind discretizations by Banks and Henshaw introduced in \cite{BANKS20125854} were generalized to overset grids. \bb{In \cite{Hagstrom2012} convergence at 11th order for a finite difference method is demonstrated.} A second order accurate overset grid method for elastic waves can be found in \cite{smog}.

We use the recently introduced dissipative Hermite methods for the scalar wave equation in second order form, \cite{secondHermite}, for the approximation on Cartesian grids. To handle geometry we use the energy based DG methods of \cite{Upwind2} on thin grids that are grown out from physical boundaries. We use projection to transfer the solutions between grids rather than interpolation.

Both the Hermite and DG methods we employ increase the order of accuracy by increasing the number of degrees of freedom on an element or cell. This has practical implications for grid generation as a single grid with minimal overlap can be used independent of order, reducing the complexity of the grid generation step. This can be important for example in problems like optimal shape design, where the boundary changes throughout the optimization. This is different from the finite difference methods where, due to the wider finite difference stencils, the overlap must grow as the order is increased.

The transfer of solutions between overset grids typically causes a perturbation to the discrete operators which, especially for hyperbolic problems, results in instabilities, see \cite{smog} for example.  These instabilities are often weak and can thus be suppressed by a small amount of artificial dissipation. There are two drawbacks of this  added dissipation, first it is often not easy to determine the suitable amount needed, i.e. big enough to suppress instabilities but small enough not to reduce the accuracy or timestep too severely. Second, in certain cases the instabilities are strong enough that the dissipation must scale with the discretization parameter (the grid size) in such a way that the order of accuracy of the overall method is reduced by one.

Similar to~\cite{ANGEL2018534}, we use a dissipative method that has naturally built--in damping that is sufficient to suppress the weak instabilities caused by the overset grids. The order of the hybrid overset grid method is the design order of the Hermite method or DG method, whichever is the smallest.


In the hybrid H--DG overset grid method the Hermite method is used on a Cartesian grid  in the interior of the domain, and the discontinuous Galerkin method on another, curvilinear grid at the boundary. The numerical solution is evolved independently on these grids for one timestep of the Hermite method. By using the Hermite method in the interior the strict timestep constraints of the DG method are relaxed by a factor that grows with the order of the method. Asymptotically, as discussed above, the complexity of the hybrid H--DG solver approaches that of the Cartesian Hermite solver~\cite{secondHermite}.

The paper is organized as follows. The Hermite method is described in the next section. We first explain the method in simple one dimensional case and then explain how the method generalizes to two dimensions. The DG method is described in section~\ref{sec:dg}. The details of the overset grids and a hybridization of the DG and the Hermite methods are described in section~\ref{sec:overset}. We illustrate the hybrid H--DG method with numerical simulations in the section~\ref{sec:experiment}.

\section{Dissipative Hermite method for the scalar wave equation} \label{sec:Hermite}
We present the Hermite method in some detail here and refer the reader to the original work~\cite{secondHermite} for convergence analysis and error estimates.

Consider the one dimensional wave equation in second order form in space and first order in time
\begin{eqnarray}
u_t &=& v,  \label{eq1a}  \\
v_t &=& c^2 u_{xx} + f, \ \  x\in \Omega, \ \  t \in (0, T). \label{eq1b}
\end{eqnarray}
Here $u \in C^{2m+3}\left(\Omega \times [0, T]\right)$, $v \in C^{2m+1} \left(\Omega \times [0,T] \right)$ and $f \in C^{2m+1}(\Omega \times [0, T])$ for optimal convergence.
We refer to $u$ as the displacement, and $v$ as the velocity. The speed of sound is $c$. We consider boundary conditions of Dirichlet or Neumann type
\bea
u(t,x) &=& h_0(t,x), \ \ x \in \partial \Omega_D, \nonumber \\
u_x(t,x) &=& h_1(t,x), \ \ x \in \partial \Omega_N, \nonumber
\eea
and initial conditions
\begin{eqnarray*}
    u(0,x) &=& g_0(x), \\
    v(0,x) &=& g_1(x).
\end{eqnarray*}

Let the spatial domain be $\Omega = [a, b]$. The domain will be discretized by a primal grid
\[
x_{i} = a + ih, \ \ h = (b-a)/N,\ \ i = 0,\dots, N,
\]
and a dual grid
\[
x_i = a + ih, \ \, i = \hf,\dots,N-\hf.
\]
The use of staggered grids \rr{allow us to evaluate the derivatives of the polynomial approximations to derivatives at the cell center, rather than throughout the cell as in most other element based methods. The slow growth with polynomial degree of the derivative approximations near the cell centers (see \cite{TAMECFL}) allows us to use timesteps that are bounded by the speed of sound and not by the degree of the polynomial.}
In time we discretize using a uniform grid with increments $\Delta t / 2$, that is
\be
t_n = n \Delta t, \ n = 0,1/2,1,\ldots. \nonumber
\ee
At each grid point $x_i$ the approximation to the solution is represented by its degrees of freedom (DOF) that approximate the values and spatial derivatives of $u$ and $v$. Equivalently, the approximations to $u$ and $v$ can be represented as polynomials centered at grid points $x_i$. The Taylor coefficients of these polynomials are scaled versions of the degrees of freedom. To achieve the optimal order of accuracy $(2m+1)$ we require the $(m+1)$ and $m$ first derivatives of $u$ and $v$ respectively to be stored at each grid point.

At the initial time (which we take to be $t = 0$) these polynomials are approximations to the initial condition on the primal grid
\begin{align}
  u(x,0) &\approx \sum_{l=0}^{m+1} \hat{u}_{l} \bigg(\frac{x-x_i}h\bigg)^l \equiv p_i(x), \ \ i = 0,\dots,N,\nonumber \\
  v(x,0) &\approx \sum_{l=0}^{m}  \hat{v}_l \bigg(\frac{x-x_i}h\bigg)^l \equiv q_i(x),
\ \ i = 0,\dots,N. \nonumber
\end{align}
The coefficients $\hat{u}_{l}$ and $\hat{v}_{l}$ are assumed to be accurate approximations to the scaled Taylor coefficients of the initial data. If expressions for the derivatives of the initial data are known we simply set
\be
\hat{u}_l = \frac{h^l}{l!}\frac{d^l g_0}{dx^l}\bigg|_{x=x_i},\ \
\hat{v}_l = \frac{h^l}{l!}\frac{d^l g_1}{dx^l}\bigg|_{x=x_i}. \label{proj}
\ee
Alternatively, if only the functions $g_0$ and $g_1$ are known, we may use a projection or interpolation procedure to find the coefficients in (\ref{proj}).

The numerical algorithm for a single timestep consists of two phases, an interpolation step and an evolution step. First, during the interpolation phase the spatial piecewise polynomials are constructed to approximate the solution at the current time. Then, in the time evolution phase we use the spatial derivatives of the interpolation polynomials to compute time derivatives of the solution using the PDE. We compute new values of the DOF on the next time level by evaluating the obtained Taylor series.  We now describe each step separately.

\subsection{Hermite interpolation}\label{ssec:Hermite:int}
At the beginning of a timestep at time $t_n$ (or at the initial time) we consider a cell $[x_i, x_{i+1}]$ and construct the unique local Hermite interpolant of degree $(2m+3)$ for the displacement and degree $(2m+1)$ for the velocity. The interpolating polynomials are centered at the dual grid points $x_{i+\hf}$ and can be written in Taylor form
\bea
p_{i+\hf}(x) &=& \sum_{l=0}^{2m+3} \hat{u}_{l,0} \bigg(\frac{x-x_{i+\hf}}h\bigg)^l,\ \ x\in[x_{i}, x_{i+1}],
\ \ i = 0,\dots,N-1, \label{hi:p} \\
q_{i+\hf}(x) &=& \sum_{l=0}^{2m+1} \hat{v}_{l,0} \bigg(\frac{x-x_{i+\hf}}h\bigg)^l,\ \ x\in[x_{i}, x_{i+1}],
\ \ i = 0,\dots,N-1. \label{hi:q}
\eea
The interpolants $p_{i+\hf}$ and $q_{i+\hf}$ are determined by the local interpolation conditions:
\begin{align}
\frac{d^l p_{i+\hf}}{dx^l} &= \frac{d^l p_i}{dx^l}\bigg|_{x=x_i},\ \
\frac{d^l p_{i+\hf}}{dx^l} = \frac{d^l p_{i+1}}{dx^l}\bigg|_{x=x_{i+1}}, \ \ l = 0,\dots,m+1,
\nonumber \\
\frac{d^l q_{i+\hf}}{dx^l} &= \frac{d^l q_i}{dx^l}\bigg|_{x=x_i},\ \
\frac{d^l q_{i+\hf}}{dx^l} = \frac{d^l q_{i+1}}{dx^l}\bigg|_{x=x_{i+1}}, \ \ l = 0,\dots,m.
\nonumber
\end{align}
We find the coefficients in \eqref{hi:p} and \eqref{hi:q} by forming a generalized Newton table as described in \cite{hagstrom2015solving}.

\subsection{Time evolution}\label{ssec:Hermite:ts}
To evolve the solution in time we further expand the coefficients of $p_{i+\hf}$ and $q_{i+\hf}$. At each point on the dual grid, $x_{i+\hf}$ we seek temporal Taylor series
\bea
p_{i+\hf}(x,t) &=& \sum_{l=0}^{2m+3}\sum_{s=0}^{\kappa_p}
\hat{u}_{l,s} \bigg(\frac{x-x_{i+\hf}}h\bigg)^l\bigg(\frac{t}{\Delta t}\bigg)^s,\label{tts:p} \\
q_{i+\hf}(x,t) &=& \sum_{l=0}^{2m+1}\sum_{s=0}^{\kappa_q}
\hat{v}_{l,s} \bigg(\frac{x-x_{i+\hf}}h\bigg)^l\bigg(\frac{t}{\Delta t}\bigg)^s,\label{tts:q}
\eea
where $\kappa_p = (2m+3-2\lceil \frac l2 \rceil)$ and $\kappa_q=(2m+1-2\lfloor \frac l2 \rfloor)$. The coefficients $ \hat{u}_{l,0}$ and $ \hat{v}_{l,0}$ are given
by the coefficients of \eqref{hi:p} and \eqref{hi:q}. At this time the scaled time derivatives, $ \hat{u}_{l,s}$ and $ \hat{v}_{l,s}$ $s > 0$, are unknown and must be determined. Once they are determined we may simply evaluate (\ref{tts:p}) and (\ref{tts:q}) at $t = t_n + \Delta t /2$ to find the solution at the next half timestep.

In Hermite methods the coefficients of temporal Taylor polynomials are determined by collocating the differential equation, \cite{good2006,secondHermite,hagstrom2015solving}. In particular, by differentiating \eqref{eq1a} and \eqref{eq1b} in space and time the time derivatives of the solution can be directly expressed in terms of spatial derivatives
\begin{eqnarray}
    \frac{\pa^{s+1+r} u }{\pa t^{s+1} \pa x^{r}} &=& \frac{\pa^{s+r} v}{\pa t^{s} \pa x^{r}},  \label{deq1a}  \\
    \frac{\pa^{s+1+r} v }{\pa t^{s+1} \pa x^{r}} &=& c^2 \frac{\pa^{s + r + 2} u }{\pa t^{s} \pa x^{r+2}}  + \frac{\pa^{s+r} f }{\pa t^{s} \pa x^{r}}. \label{deq1b}
\end{eqnarray}
Substituting \eqref{tts:p} and \eqref{tts:q} into \eqref{deq1a}
and \eqref{deq1b} and evaluating at $x = x_{i+\hf}$ and $t = t_n$,
we can match the powers of the coefficients to find the recursion relations
\begin{eqnarray}
    \hat{u}_{l,s+1}  &=& \frac{\Delta t}{s} \hat{v}_{l,s},  \label{tcofsa}  \\
    \hat{v}_{l,s+1} &=& c^2\frac{(l+1)(l+2)}{h^2} \frac{\Delta t}{s} \hat{u}_{l+2,s} +  \frac{\Delta t}{s} \hat{f}_{l,s}. \label{tcofsb}
\end{eqnarray}
Here $\hat{f}_{l,s}$ are the coefficients of the Taylor expansion of $f$, or of the polynomial which interpolates $f(t, x_{i+1/2})$ in time around $t = t_n$. Note that since there are a finite number of coefficients, representing the spatial derivatives at the time $t_n$, the recursions truncate and only $\kappa_p$ and $\kappa_q$ terms need to be considered.

To complete a half timestep we evaluate the approximation at $t =t_n+ \frac{\Delta t}{2}$ for the $(m+1)$ and $m$ first derivatives
\begin{align}
& \frac{\pa^l u}{\pa x^l} (x_{i+\hf},t_{n+\frac{1}{2}})  \approx \frac{\pa^l p_{i+\hf}}{\pa x^l}  (x_{i+\hf},t_{n+\frac{1}{2}})
= \frac{l!}{h^l} \sum_{s=0}^{\kappa_p} \frac{\hat{u}_{l,s}}{2^s} , \ \ l = 0,\dots, m+1,
\label{tsa}\\
& \frac{\pa^l v}{\pa x^l} (x_{i+\hf},t_{n+\frac{1}{2}})  \approx \frac{\pa^l q_{i+\hf}}{\pa x^l} (x_{i+\hf},t_{n+\frac{1}{2}})
=\frac{l!}{h^l}
 \sum_{s=0}^{\kappa_q} \frac{\hat{v}_{l,s}}{2^s}, \ \ l = 0,\dots, m. \label{tsb}
\end{align}

\begin{rem}
\rr{A remarkable feature of Hermite methods is that (independent of order of accuracy) since the initial data for each cell is a polynomial the time evolution is exact whenever the following conditions are met: 1.) The recursion relations (\ref{tcofsa}) and (\ref{tcofsb}) are run until they truncate, 2.) The forcing is zero (or a polynomial of degree $2m+1$), 3.) Each cell $[x_i, x_{i+1}]$ includes the base of the domain of  dependence of the solution at a dual grid point $x_{i+\hf}$ at time $t = \frac {\Delta t}2$ (see e.g. \cite{secondHermite}). The latter condition can also be stated as a CFL condition}
\[
c\frac{\Delta t}{2} \le \frac{h}{2}.
\]
\rr{
In the present method we do not quite achieve this optimal CFL condition but have verified numerically that our solvers of orders of accuracy 3, 5 and 7 are stable for $c\Delta t \le 0.75 h$.}
\end{rem}

\bb{
\subsubsection{Variable coefficients}
For problems with a variable wave speed the acceleration is governed by
\begin{equation}
  v_t = \underbrace{(c^2(x) u_x)_x}_{s(x)}. \label{eq:vc1}
\end{equation}
To compute $v_t, v_{tt}, v_{ttt}, $ etc., needed to evolve the solution by a Taylor series method, we must evaluate the right hand side of (\ref{eq:vc1}). This is done in sequence by forming the polynomial $s(x)$ by operations on polynomials. Let $p \approx u$ and $a \approx c^2$ be polynomials approximating $u(t,x)$ and $c^2(x)$ then
\[
s(x) = \mathcal{D} (a(x) \otimes (\mathcal{D} p(t,x))).
\]
Here $\mathcal{D}$ denotes polynomial differentiation and $\otimes$ represents polynomial multiplication with degree truncation to the degree of $p$. Computation of $v_t, v_{tt}$ etc. can now be done by forming $s(x)$ with $p \approx u, u_t, u_{tt}$, etc.
}

\subsection{Imposing boundary conditions for the Hermite method}\label{ssec:Hermite:bc}
\bb{In the hybrid Hermite-DG overset grid method, physical boundary conditions can be imposed on any grid that discretizes the boundary. For example, in the numerical experiments in Section \ref{ss:body}, the boundary conditions are imposed on both grids. In this section we explain how physical boundary conditions are imposed for the Hermite method and a Cartesian grid.}

Physical boundary conditions are enforced at the half time level, i.e. when the solution on the dual grid is to be advanced back to the primal grid. As there are many degrees of freedom that are located on the boundary the physical boundary condition must be augmented by the differential equation to generate more independent conditions so that the degrees of freedom can be uniquely determined. The basic principle, often referred to as compatibility boundary conditions (see e.g.  \cite{henshaw:1730}), is to take tangential derivatives of the boundary conditions and combine these with the PDE.

For example, assume we want to impose the boundary condition
\be
u(t,0) = g(t). \label{e:bc}
\ee
\rr{
Then, as $x_0=0$ is a boundary grid point the Taylor polynomials
}
\eqref{tts:p}-\eqref{tts:q},
\rr{centered at $x_0$, should satisfy the boundary condition
}
\eqref{e:bc}
\rr
{and compatibility conditions, (i.e. conditions for the derivatives), that one obtains by differentiating} \eqref{e:bc} \rr{in time and then replace time derivatives of $u$ in favor of spatial derivatives by using the wave equation}.
We thus seek a polynomial outside the domain which together with the polynomial just inside the boundary forms a Hermite interpolant that satisfies the boundary and compatibility conditions.

\rr{
Precisely, to evolve the solution on the boundary we must determine the $2(m+2)$ and $2(m+1)$ coefficients of the polynomials approximating $u$ and $v$ at the boundary. For example for $u$, this polynomial must interpolate the $(m+2)$ data describing the current approximation of $u$ at dual grid point next to the boundary, this yields $(m+2)$ independent linear equations. The remaining $(m+2)$ independent linear equations can be obtained by requiring that the polynomial satisfies with the boundary condition $u(0,t) =  g(t)$ and its time derivatives as described above.
}

Once the interpolant is determined on the boundary we evolve it as in the interior (see section \ref{ssec:Hermite:ts}).

\begin{rem}
  We note that in the special case of a flat boundary and homogeneous Dirichlet or Neumann boundary conditions then enforcing the boundary conditions reduces to enforcing that the polynomial on the boundary is either odd or even, respectively, in the normal direction. Then the correct odd polynomial can be obtained by constructing the polynomial outside the domain $\Omega$ (often referred as ghost-polynomial) by mirroring the coefficients corresponding to even powers in the normal coordinate variable with a negative sign and the  coefficients corresponding to odd powers with the same sign.
\end{rem}

Boundary conditions at interior overset grid boundaries are supplied by projection of the known solutions from  other grids and will be discussed below.

\subsection{Higher dimensions}
In higher dimensions the approximations to $u$ and $v$ take the form of centered tensor product Taylor polynomials. In two dimensions (plus time) the coefficients would be of the form \bb{$\hat u_{k,l,s}$}, with the two first indices representing the powers in the two spatial directions, and the third representing time.

For the scalar wave equation
\begin{align*}
u_t &= v, \\
v_t &= c^2 (u_{xx}+u_{yy}), \ \  (x,y) \in \Omega, \ \  t > 0,
\end{align*}
the recursion relations for computing the time derivatives are a straightforward generalization of the one dimensional case
 \begin{align}
& \bb{\hat u}_{k,l,s} = \f{\Delta t}{s} \bb{\hat v}_{k, l, s-1}, \label{eq:recursion_ut} \\
& \bb{\hat v}_{k, l, s} =
c^2\f{(k+2)(k+1)}{s} \f{\Delta t}{h_x^2} \bb{\hat u}_{k+2, l, s-1} +
c^2\f{(l+2)(l+1)}{s} \f{\Delta t}{h_y^2} \bb{\hat u}_{k, l+2,  s-1}. \label{eq:recursion_vt}
\end{align}
As noted in \cite{secondHermite}, using this recursion for all the time derivatives does not produce a method with order independent CFL condition but a method whose time-step size decrease slightly as the order increases. For optimally large timesteps it is necessary to use the special start up procedure
\begin{align*}
  \bb{\hat u}_{k,l,1} &= \Delta t\, \bb{\hat v}_{k, l,0}, \\
  \bb{\hat v}_{k, l, 1} &=  \Delta t c^2 \left( \frac{(k+2)(k+1)}{h_x^2} \bb{\hat u}^X_{k+2, l} +
  \frac{(l+2)(l+1)}{h_y^2} \bb{\hat u}^Y_{k,l+2} \right).
\end{align*}
Here $\bb{\hat u}^X_{k, l}$ are the $(2m+4) \times (2m+2)$ coefficients of the interpolating polynomial of degree $(2m+3)$ in $x$ and degree $(2m+1)$ in $y$ and $\bb{\hat u}^Y_{k, l}$ are the $(2m+4) \times (2m+2)$ coefficients of the interpolating polynomial of degree $(2m+3)$ in $y$ and degree $(2m+1)$ in $x$. For the remaining coefficients $s = 2\ldots, 4m+3$ we use (\ref{eq:recursion_ut}) and (\ref{eq:recursion_vt}) with $k,l = 0,\ldots,2m+1$. Further details of the two dimensional method can be found in \cite{secondHermite}.

\section{Energy based discontinuous Galerkin methods for the wave equation} \label{sec:dg}
Our spatial discontinuous Galerkin discretization is a direct application of the energy based formulation described for general second order wave equations in \cite{Upwind2,el_dg_dath,fluid_solid_DG}. Here, our energy based DG method starts from the energy of the scalar wave equation
\[
H(t) = \int_{\Omega} \frac{v^2}{2} + G(x,y,\nabla u) \, d \Omega,
\]
where $ G(x,y,\nabla u) = \frac{c^2(x,y)}{2} |\nabla u|^2 $ is the potential energy density, $v$ is the velocity or the time derivative of the displacement, $v = u_t$.

 Now, the wave equation, written as a second order equation in space and first order in time takes the form
\[
u_t = v, \ \ v_t = - \delta G,
\]
where $\delta G$ is the variational derivative of the potential energy
\[
\delta G = - \nabla \cdot (c^2(x,y) \nabla u).
\]
For the continuous problem the change in energy is
\be
\f{d H(t)}{dt} = \int_{\Omega} v v_t + u_t  \left[ \nabla \cdot (c^2(x,y) \nabla u \right] \,d \Omega = [ u_t  (n \cdot (c^2(x,y) \nabla u))]_{\partial \Omega}, \nonumber
\ee
where the last equality follows from integration by parts together with the wave equation.

A variational formulation that mimics the above energy identity can be obtained if the equation $v-u_t=0$ is tested with the variational derivative of the potential energy. Let $\Omega_j$ be an element and $(\Pi^{q_u}(\Omega_j))^2$ and $(\Pi^{q_v}(\Omega_j))^2$ be the spaces of tensor product polynomials of degrees $q_u$ and $q_v=q_u-1$. Then, the variational formulation on that element is:
\begin{problem}
  Find $v^h \in (\Pi^{q_v}(\Omega_j))^2$, $u^h \in (\Pi^{q_u}(\Omega_j))^2$ such that for all
  $\psi \in (\Pi^{q_v}(\Omega_j))^2$, $\phi \in (\Pi^{q_u}(\Omega_j))^2$
\begin{align}
\int_{\Omega_j} (c^2 \nabla \phi) \cdot \left( \f {\pa \nabla u^h}{\pa t}-\nabla v^h \right) d \Omega & =
[ (c^2 \nabla \phi ) \cdot n \left( v^{\ast}-v^h \right)]_{\pa \Omega_j},
\label{var1} \\
\int_{\Omega_j} \psi \f {\pa v^h}{\pa t} + c^2 \nabla \psi \cdot  \nabla u^h \, d\Omega& =
 [\psi \, (c^2\,\nabla u \cdot n)^{\ast}]_{\pa \Omega_j}. \label{var2}
\end{align}
\end{problem}

Let $[[f]]$ and $\{f\}$ denote the jump and average of a quantity $f$ at the interface between two elements, then, choosing the numerical fluxes as
\begin{align}
v^{\ast}  &= \{v^h\} -\tau_1 [[c^2\,\nabla u^h \cdot n]], \nonumber\\
(c^2\,\nabla u \cdot n)^{\ast} &= \{ c^2\,\nabla u^h \cdot n \}  -\tau_2 [[v^h]], \nonumber
\end{align}
yields a contribution $ -\tau_1 ([[c^2\,\nabla u^h \cdot n]])^2 -\tau_2 ([[v^h]])^2$ from each element face to the change of the discrete energy, guaranteeing that
\[
\f{d H^h(t)}{dt} \equiv \frac{d}{dt} \sum_{j} \int_{\Omega_j} \frac{(v^h)^2}{2} + G(x,y,\nabla u^h) \le 0.
\]
Physical boundary conditions are enforced through the numerical fluxes, see \cite{Upwind2} for details.

Note that the above energy estimate follows directly from the formulation \eqref{var1} - \eqref{var2} but as the energy is invariant to constants equation \eqref{var1} must be supplemented by the equation
\begin{equation}
\int_{\Omega_j} \left( \f {\pa u^h}{\pa t} - v^h \right) d \Omega  = 0. \nonumber 
\end{equation}

Our implementation uses quadrilaterals and approximations by tensor product Chebyshev polynomials of the solution on the reference element $(r,s) \in [-1,1]^2$. That is, on each quadrilateral we have approximations on the form
\begin{align}
  u(x(r,s),y(r,s),t_n) &\approx  \sum_{l=0}^{q_u}
\sum_{k=0}^{q_u} c_{lk}
T_l (r)
T_k(s), \nonumber \\
  v(x(r,s),y(r,s),t_n) &\approx  \sum_{l=0}^{q_v}
\sum_{k=0}^{q_v} d_{lk}
T_l (r)
T_k(s). \nonumber
\end{align}
We choose $\tau_1 = \tau_2 = 1/2$ (so called upwind or Sommerfeld fluxes) which result in methods where $u$ is observed to be $q_u+1$ order accurate in space \cite{Upwind2}. \bb{We note that another basis like Legendre polynomials could also be used. In fact we have repeated some of the long time  computations in the numerical experiments section below to confirm that a change of basis to Legendre polynomials does not effect the stability or accuracy properties of the method.}

\subsection{Taylor series time-stepping}\label{ssec:DG:ts}
In order to match the order of accuracy in space and time for the DG method we employ Taylor series time-stepping.  Assuming that all the degrees of freedom have been assembled into a vector ${\bf w}$ we can write the semi-discrete method as ${\bf w}_t = A {\bf w} $ with $A$ being the matrix  representing the spatial discretization. If we know the discrete solution at the time $t_n$ we can advance it to the next time step $t_{n+1} = t_n + \Delta t$ by the simple formula
\begin{eqnarray*}
{\bf w}(t_{n}+\Delta t) &=& {\bf w}(t_{n}) + \Delta t {\bf w}_t(t_{n}) +  \frac{(\Delta t)^2}{2!}{\bf w}_{tt}(t_{n}) \ldots \\
 &=& {\bf w}(t_{n}) + \Delta t A {\bf w}(t_{n}) +  \frac{(\Delta t)^2}{2!} A^2 {\bf w}(t_{n}) \ldots
\end{eqnarray*}
As we use dissipative fluxes this timestepping method is stable as long as the number of stages in the Taylor series is greater than the order of accuracy in space and with the timestep small enough.

\section{Overset grid methods} \label{sec:overset}
In this section we explain how we use the two discretization techniques described above on overset grids to approximate solutions to the scalar wave equation.

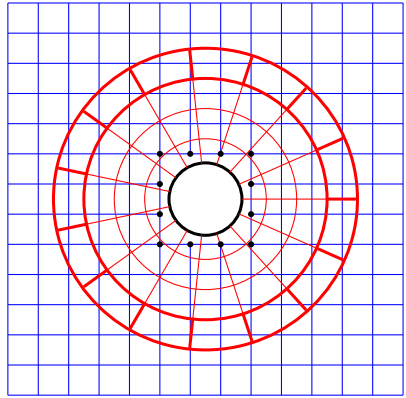
\begin{figure}[htb]
\begin{center}
 \setlength{\unitlength}{0.4mm}
  \begin{picture}(130,130)(-65,-65)
    \multiput(-65,-65)(10,0){6}{\color{blue}\line(0,1){130}}
    \multiput(-65,-65)(0,10){6}{\color{blue}\line(1,0){130}}
    \multiput(65,65)(-10,0){6} {\color{blue}\line(0,-1){130}}
    \multiput(65,65)(0,-10){6} {\color{blue}\line(-1,0){130}}
    \multiput(5,15)(-10,0){2}  {\color{blue}\line(0,1){50}}
    \multiput(5,-15)(-10,0){2} {\color{blue}\line(0,-1){50}}
    \multiput(15,5)(0,-10){2}  {\color{blue}\line(1,0){50}}
    \multiput(-15,5)(0,-10){2} {\color{blue}\line(-1,0){50}}
    \put(0,0){\color{black}\linethickness{0.4mm}\circle{24}}
    \put(0,0){\color{red}\circle{40}}
    \put(0,0){\color{red}\circle{60}}
    \put(0,0){\color{red}\linethickness{0.4mm}\circle{80}}
    \put(0,0){\color{red}\linethickness{0.4mm}\circle{100}}
      \put( 12.0000,        0){\color{red}\line( 12.0000,        0){38}}
      \put( 10.9625,   4.8808){\color{red}\line( 10.9625,   4.8808){34.5}}
      \put(  8.0296,   8.9177){\color{red}\line(  8.0296,   8.9177){25.5}}
      \put(  3.7082,  11.4127){\color{red}\line(  3.7082,  11.4127){11.7}}
      \put( -1.2543,  11.9343){\color{red}\line( -1.2543,  11.9343){4}}
      \put( -6.0000,  10.3923){\color{red}\line( -6.0000,  10.3923){19}}
      \put( -9.7082,   7.0534){\color{red}\line( -9.7082,   7.0534){30.6}}
      \put(-11.7378,   2.4949){\color{red}\line(-11.7378,   2.4949){37.5}}
      \put(-11.7378,  -2.4949){\color{red}\line(-11.7378,  -2.4949){37.5}}
      \put( -9.7082,  -7.0534){\color{red}\line( -9.7082,  -7.0534){30.6}}
      \put( -6.0000, -10.3923){\color{red}\line( -6.0000, -10.3923){19}}
      \put( -1.2543, -11.9343){\color{red}\line( -1.2543, -11.9343){4}}
      \put(  3.7082, -11.4127){\color{red}\line(  3.7082, -11.4127){11.7}}
      \put(  8.0296,  -8.9177){\color{red}\line(  8.0296,  -8.9177){25.5}}
      \put( 10.9625,  -4.8808){\color{red}\line( 10.9625,  -4.8808){34.5}}

      \put(   40.0000,        0){\linethickness{0.4mm}\color{red}\line( 12.0000,        0){ 9.5000  }}
      \put(   36.5418,  16.2695){\linethickness{0.4mm}\color{red}\line( 10.9625,   4.8808){   8.6250}}
      \put(   26.7652,  29.7258){\linethickness{0.4mm}\color{red}\line(  8.0296,   8.9177){   6.3750}}
      \put(   12.3607,  38.0423){\linethickness{0.4mm}\color{red}\line(  3.7082,  11.4127){   2.9250}}
      \put(   -4.1811,  39.7809){\linethickness{0.4mm}\color{red}\line( -1.2543,  11.9343){   1.0000}}
      \put(  -20.0000,  34.6410){\linethickness{0.4mm}\color{red}\line( -6.0000,  10.3923){   4.7500}}
      \put(  -32.3607,  23.5114){\linethickness{0.4mm}\color{red}\line( -9.7082,   7.0534){   7.6500}}
      \put(  -39.1259,   8.3165){\linethickness{0.4mm}\color{red}\line(-11.7378,   2.4949){   9.3750}}
      \put(  -39.1259,  -8.3165){\linethickness{0.4mm}\color{red}\line(-11.7378,  -2.4949){   9.3750}}
      \put(  -32.3607, -23.5114){\linethickness{0.4mm}\color{red}\line( -9.7082,  -7.0534){   7.6500}}
      \put(  -20.0000, -34.6410){\linethickness{0.4mm}\color{red}\line( -6.0000, -10.3923){   4.7500}}
      \put(   -4.1811, -39.7809){\linethickness{0.4mm}\color{red}\line( -1.2543, -11.9343){   1.0000}}
      \put(   12.3607, -38.0423){\linethickness{0.4mm}\color{red}\line(  3.7082, -11.4127){   2.9250}}
      \put(   26.7652, -29.7258){\linethickness{0.4mm}\color{red}\line(  8.0296,  -8.9177){   6.3750}}
      \put(   36.5418, -16.2695){\linethickness{0.4mm}\color{red}\line( 10.9625,  -4.8808){   8.6250}}

      \multiput(-15,-15)(10,0){4}{\circle*{2}}
      \multiput(-15,15)(10,0){4}{\circle*{2}}
      \multiput(-15,-5)(0,10){2}{\circle*{2}}
      \multiput(15,-5)(0,10){2}{\circle*{2}}

  \end{picture}
  \caption{An example of an overset grid for a circular boundary inside a square. The red grid is curvilinear and the blue grid is Cartesian (in a realistic problem the red grid would be significantly thinner). The black filled circles indicate the cut out domain boundary. \label{fig:composite_grid}}
\end{center}
\end{figure}

The idea behind the overset grid methods is to cover the bulk of the domain with a Cartesian grid, where efficient methods can be employed, and to discretize the geometry with narrow body-fitted grids. In Figure \ref{fig:composite_grid} we display two overset grids, a blue Cartesian grid, which we denote $a$, and a red curvilinear grid, which we denote $b$, that are used to discretize a geometry consisting of a circular hole cut out from a square region. Note that the grids overlap, hence the name overset grids.  Also, note that the annular grid cuts out a part of the Cartesian grid. This cut of the Cartesian grid creates \rr{an} internal, non-physical boundary in the blue grid.

Here physical boundary conditions are enforced on the red grid at the black boundary which defines the inner circle and on the outermost boundary on the blue grid.

In order to use the Hermite or DG methods on the grids we will need to supply boundary conditions at the interior boundaries. In the example in Figure \ref{fig:composite_grid} this means that we would have to specify the solution on the outer part of the annular grid and on the staircase boundary (marked with filled black circles) that has been cut out from the Cartesian grid.

\begin{figure}[H]
  \includegraphics[width=0.49\textwidth]{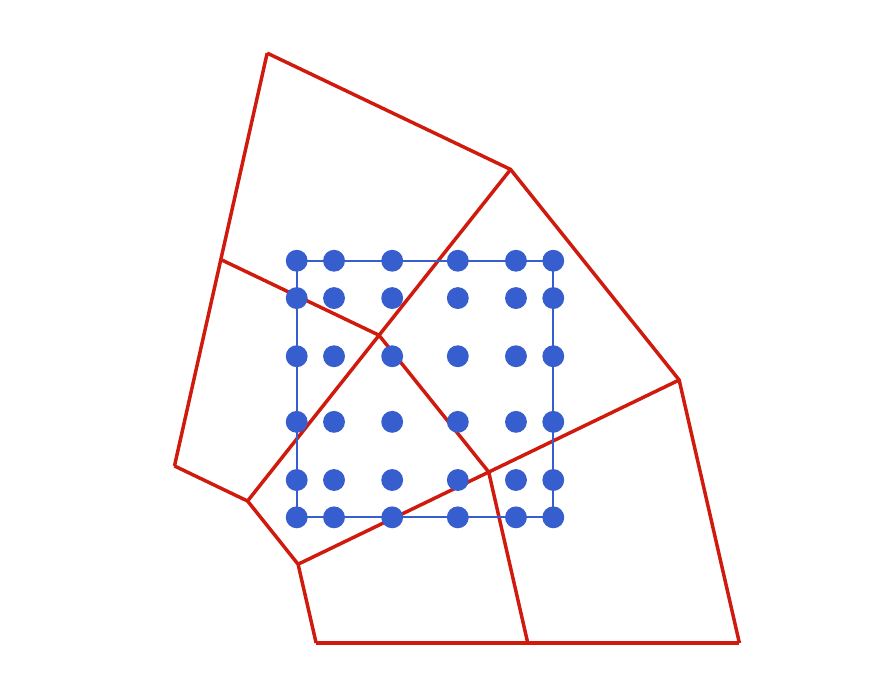}
  \includegraphics[width=0.49\textwidth]{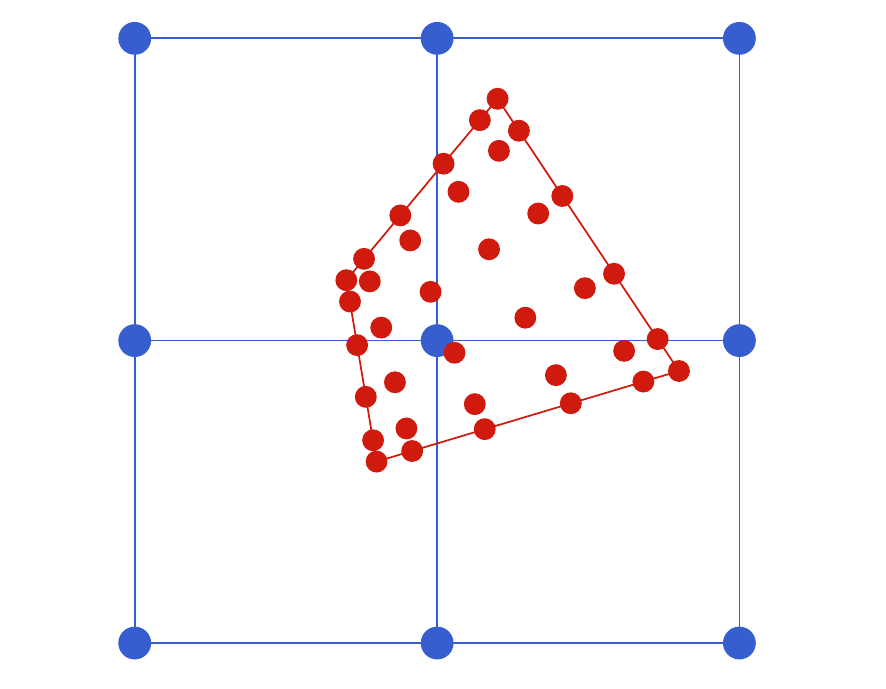}
  \caption{Typical setup for communication. In the left subfigure the local tensor product GLL grid around a Hermite grid point is marked with filled blue circles. The points in the GLL grid may be covered by different DG elements. In the right subfigure the tensor product grid inside the DG element is marked with filled red circles. The points in the GLL grid may be contained in different Hermite cells.  \label{f:ogm1}}
\end{figure}

In most methods that use overset grids, in particular those using finite differences, the communication of the solution on the interior boundaries is done by interpolation, see e.g. \cite{chess1990}. For the methods we use here we have found that the stability properties are greatly enhanced if we instead transfer volumetric data (numerical solution) in the elements / gridpoints near the internal boundaries by projection rather than by interpolation. In fact, when we use volume data the resulting methods are stable without adding artificial dissipation, when we use interpolation they are not. \rr{At the end of this section we discuss a possible reason why the projection behaves better than interpolation.}

As mentioned above, in a Hermite method, we can think of the degrees of freedom as  either being nodal data, consisting of function and derivative values, or as coefficients in a Taylor polynomial. Thus, when transferring data to a grid where a Hermite method is used (like the example in the left subfigure of Figure \ref{f:ogm1}) we must determine a tensor product polynomial centered around a gridpoint local to that grid (the points we would center around are indicated by black points in Figure \ref{fig:composite_grid}). Below we will explain in detail how we determine this polynomial.

For elements with an internal boundary face (denoted by thick red lines in Figure~\ref{fig:composite_grid}) we could in principle transfer the solution by specifying a numerical flux on that face, however we have found that this approach results in weakly unstable methods. Instead we transfer volumetric data to each element that has an internal boundary face, we give details below. Given the timestep constraints of DG methods we must march the DG solution using much smaller timesteps than those used for the Hermite method. This necessitates the evaluation of the Hermite data not only at the beginning of a Hermite timestep but at many intermediate times.

\subsection{Determining internal boundary data for the Hermite solver}\label{ssec:Overset:hdg}
We first consider the problem of determining internal boundary data required by the Hermite method. An example of how to compute solution data at the gridpoints $(x_i, y_{j})$ at the boundary of Cartesian grid (filled black circles) is depicted in Figure \ref{fig:composite_grid}.

In general, the tensor product polynomial centered around $(x_i, y_{j})$
is found by a two step procedure. First we project into a local $L_2$ basis spanned by Legendre polynomials and perform a numerically stable and fast change of basis into the monomial basis. Then we truncate the monomial to the degree required by the Hermite method.

To carry out the $L_2$ projection we introduce a local tensor product Gauss-Legendre-Lobatto (GLL) grid centered around $(x_{i}, y_{j})$. These points are marked as filled blue circles in the left subfigure of Figure~\ref{f:ogm1}. The number of grid points in the local grids are determined by the order of the projection. To maintain the order of the method, the order of the projection should be at least the same as the order of the spatial discretization, thus it is sufficient to have $2m+4$ points in each direction. The GLL quadrature nodes are defined on the reference element $(r,s) \in [-1,1]^2$ that maps to a cell defined by the dual gridpoints closest to $(x_{i}, y_{j})$.

Let $\tilde{u}$ be the numerical solution on the red grid. In the first step of the communication we compute the coefficients of a polynomial $\tilde{p}$ approximating $\tilde{u}$ by projecting $\tilde{u}$ on the space of tensor product Legandre polynomials $P_lP_k$, that is
\be
\tilde{p}(r,s) = \sum_{l=0}^{2m+3}
\sum_{k=0}^{2m+3} c_{lk}
P_l (r)
P_k(s),\
c_{lk} =
\frac{(\tilde{u}, P_l P_k)}
{\| P_l  P_k \|^2}. \label{eq:leg_proj}
\ee
Here $(f,g)$ denotes the $L_2$ inner product on $(r,s) \in [-1,1]^2$ and $\| f \|^2_2 = (f,f)$ is the norm induced by the inner product. Note that the expression (\ref{eq:leg_proj}) is particularly simple since the Legendre polynomials are orthogonal on the domain of integration. To do this we evaluate  $\tilde{u}$ at the underlying blue quadrature points in the left subfigure of Figure~\ref{f:ogm1}.

Once the polynomial (\ref{eq:leg_proj}) has been found we perform a change of basis into the local monomial used by the Hermite method. Such a change of basis can be done by the fast Vandermonde techniques by Bj\"{o}rk and Pereyra, see e.g. \cite{bp70,Dahlquist:2008fu}. At this stage the polynomial is of total degree $2m+3$ so the final step is to truncate it to total degree $m$ or $m+1$ depending on whether we are considering the displacement or the velocity. With the $(m+1)^2$ and $(m+2)^2$ degrees of freedom determined everywhere on a Hermite grid we may evolve the solution as described in Section \ref{sec:Hermite}.

\subsection{Determining data for DG elements with internal boundary faces}\label{ssec:Overset:dgh}
We now consider the problem of determining the data required by the DG method. Here we show how to obtain the data at a single DG element with at least one internal boundary face. As the timesteps of the DG method are significantly smaller than for the Hermite method we must repeat the transfer of data many times. We must also explicitly transfer time derivative data in order to use a Taylor series timestepping approach.

The tensor product polynomials in our implementation of the DG method are composed by the product of Chebyshev polynomials $T_j(z) = \cos (j \cos^{-1}(z))$ that are expressed on the reference element $(r,s) \in [-1,1]^2$. Precisely we seek
\be
{p}(r,s) = \sum_{l=0}^{q}
\sum_{k=0}^{q} c_{lk}
T_l (r)
T_k(s). \nonumber 
\ee
To determine such polynomials we perform a projection of the solution $u$, i.e the solution on Cartesian grid,
\begin{equation}
  c_{lk} =
  \frac{(\tilde{u}, T_l T_k)_{\rm C}}
  {\| T_l  T_k \|_{\rm C}^2},\nonumber
\end{equation}
but in this case the weighted inner product is
\begin{equation*}
(f,g)_{\rm C} = \int_{-1}^1\int_{-1}^{1}  \frac{f(r,s) g(r,s)}{\sqrt {1-r^2} \sqrt{1-s^2}}\, dr ds,
\end{equation*}
where the Chebyshev polynomials are orthogonal. To carry out this projection we use a local tensor product Chebyshev quadrature nodes, $2m+2$ in each dimension, as shown in right subfigure of Figure~\ref{f:ogm1}.

Denoting The local time levels used by the DG solver $n$th Hermite timestep are defined to be
\be
t_{n,\nu} = t_{n,0}+\nu \Delta t_b,\ \nu = 0,\dots N_{\rm DG}, \nonumber
\ee
where $\Delta t_b$, and similarly $\Delta t_a$, are timesteps taken on grids $b$ (curvilinear) and $a$ (Cartesian) respectively. For simplicity the starting local time level and the final local time level are equal to consequent timesteps on the Hermite grid, $t_n$ and $t_{n+1}$
\be
t_{n,0} = t_n,\ \ t_{n, N_{\rm DG}} = t_{n+1}. \nonumber
\ee
To transfer the solution values and the time derivatives needed at each of the quadrature points and at each $t_{n,\nu}$ we carry our the following ``start up'' procedure at $t_{n,0}$. For each of the quadrature points we re-center the Hermite interpolants closest to it and compute the time derivatives precisely by the recursion relations described in section~\ref{sec:Hermite}. We note that this is an inexpensive computation as the interpolants have already been found as a step in the evolution of the Hermite solution, the only added operation is the re-centering.

\rr{\subsection{Discussion of projection and interpolation}
One of the differences in the present method and a finite difference method is that during the transfer of data to the Hermite method there is a degree truncation of (in one dimension and for $u$) a polynomial of degree $(2m+3)$ to a polynomial of degree $(m+1)$. It is natural to ask how the truncated polynomial depends on whether projection or interpolation was used to find the un-truncated polynomial.

Suppose the same $(2m+4)$ data has been used to determine two polynomials
\[
p_{\rm interp.}(z)= \sum_{l = 0}^{2m+3} a_l z^l, \ \ z \in [-1/2,1/2],
\]
and
\[
p_{\rm proj.}(z) = \sum_{l = 0}^{2m+3} \tilde{b}_l P_l(2z) = \sum_{l = 0}^{2m+3} b_l z^l, \ \ z \in [-1/2,1/2].
\]
Then due to the orthogonality of the projected polynomial it is clear that the truncated polynomial satisfies
\[
\int_{-\frac{1}{2}}^{\frac{1}{2}} \left( \sum_{l = 0}^{m+1} \tilde{b}_l P_l(2z) \right)^2 dx \le \int_{-\frac{1}{2}}^{\frac{1}{2}} \left(p_{\rm proj.}(z) \right)^2 dx.
\]
The polynomial determined by interpolation does not satisfy a similar inequality. In fact the truncation can cause a significant increase in the $L_2$-energy. To investigate this we find
\[
 {\{a_0^\ast,\ldots,a_{2m+3}^\ast\}} ={\rm argmax}  \frac{\int_{-\frac{1}{2}}^{\frac{1}{2}} \left( \sum_{l = 0}^{m+1} a_l z^l \right)^2 dx}{\int_{-\frac{1}{2}}^{\frac{1}{2}} \left( \sum_{l = 0}^{2m+3} a_l z^l \right)^2 dx},
\]
for 10000 randomly selected initial data and for $m = 1,2,3$. The largest ratio between square of the $L_2$-norms of the truncated and un-truncated polynomials were 19, 657 and 3555 for $m=1$, $m=2$ and $m=3$ respectively. While this does not conclusively rule out that it could be possible to use interpolation it does indicate that a projection based approach is to be preferred. We stress that the cause of the problem is the combination of the truncation and interpolation and that there is therefore not obvious that there is any advantage to use projection rather than interpolation for methods that does not have truncation (like finite difference methods).
}

\section{Numerical experiments} \label{sec:experiment}
The hybrid H--DG method is empirically stable and accurate, and here we demonstrate it with numerical experiments. To test the stability of the method in one dimension we first define the amplification matrix and compute its spectral radius. To test the stability in two dimensions, where the amplification matrix will take too long to compute, we provide the long time simulation and estimate the error growth for multiple refinements. Convergence tests in one and two dimensions are done for the domains where the exact solution is known. In the second half of this section we apply the method to the domain with complex curvilinear boundary in the experiment with wave scattering of the smooth pentagonal object. Finally, in the end of this section we apply the method to the inverse problem of locating the underground cavities as the forward solver.

\subsection{Numerical stability test}
Unlike the Hermite and DG methods, stability of the hybrid H--DG method cannot easily be shown analytically. \bb{As a weaker alternative, the stability can be investigated numerically by looking at the spectrum of the amplification matrix associated with the method, \cite{Strik}.

To construct the amplification matrix we apply the method to initial data composed of the unit vectors. The vector that is returned after one timestep is then placed as columns in a square matrix. If the spectral radius of the amplification matrix is smaller than $1$, or if the eigenvalues with magnitude one correspond to no-trivial Jordan blocks, then the amplification matrix is power-bounded.}

We consider the wave equation \eqref{eq1a}-\eqref{eq1b} on the unit interval $x \in [0,1]$ with homogeneous Dirichlet and Neumann boundary conditions
at $x = 0$ and $x = 1$ respectively.
We introduce two uniform Cartesian grids which overlap inside a small interval close to one of the boundaries.
Precisely, the grids are
\begin{align}
\Omega_a &= \{x^a_i = ih_a,\ \ i = 0, \dots, n_a\},\nonumber \\
\Omega_b &= \{x^b_i = 1 - (n_b-i)h_b, \ \ i = 0, \dots, n_b\}. \nonumber
\end{align}
The Hermite method is used on a grid $a$ and the DG method is used on grid $b$.
The grids thus overlap inside the interval $[x^b_0, x^a_{n_a}]$. Here the ratio of the overlap size and the discretization width is $(x^a_{n_a} - x^b_0)/h_a$. This ratio is fixed for all values of $h_a$ and $h_b$.  We also fix $n_b$ so that the amount of work
done on grid $b$ is constant per timestep for all refinements. Fixing the ratio $(x^a_{n_a} - x^b_0)/h_a$ and $n_b$  makes the efficiency of the overall method
to asymptotically be the determined by the efficiency the Hermite method.

Let ${\mathbf w}^n$ be a vector holding the degrees of freedom of both methods at $n$th timestep, then we may express the complete timestep evolution as ${\mathbf w}^{n+1} = {\cal H} \mathbf w^n$ where $\cal H$ incorporates timestepping and projection. $\mathcal{H}$ can be expressed as the matrix $H$ that can be computed column by column via
\be
H_k = \mathcal{H} e_k,\label{constr}
\ee
where $e_k$ is the $k$th unit vector. The equation
\be
{\mathbf w}^n = H^n {\mathbf w}^0,
\label{e:5.1.1}
\ee
is equivalent to the $n$ timesteps of the hybrid H--DG method. \bb{Let, $\lambda = \rho(H)$ is the spectral radius of $H$. If $|\lambda|<1$ then $\| H^n \|_2$ will tend to zero for large $n$. Of course this only means that this particular discretization of this particular problem is stable and does (in principle) not tell us anything about other grid configurations.}


We consider the case $c=1$ and take the parameters to be
\[
    n_a = 10, 20, \dots, 60,\ \ n_b = 5, \ \ \frac{h_b}{h_a} = 0.9.
\]
Other parameters are $q_u,q_v$ for the DG method and $n_{DG}$, the number of timesteps done by the DG method during one step of the Hermite method. The parameters $q_u$ and $q_v$ are set so the methods used have the same order of accuracy as the approximation of $v$ for the Hermite method
\[
  q_u = 2m+2,\ q_v = 2m+1.
 \]
To get an optimal $n_{DG}$, we take the largest possible timestep for the energy based DG method \rr{(empirically determined in \cite{Upwind2})}, so that
\[
  \f{\Delta t_b}{h_b} \le 0.15/q_u,
\]
and
\[
n_{\rm DG} = \f {\Delta t_a}{\Delta t_b},
\]
is an integer. Equivalently, if the Hermite method CFL number is set, we get
\[
  n_{\rm DG} = \frac{\Delta{t_a}}{\Delta{t_{b}}} = \left\lceil {\rm CFL }\frac{q_u}{0.15}\f{h_a}{h_b} \right\rceil.
\]

\begin{figure}[]
    \begin{center}
    \includegraphics[width=\textwidth]{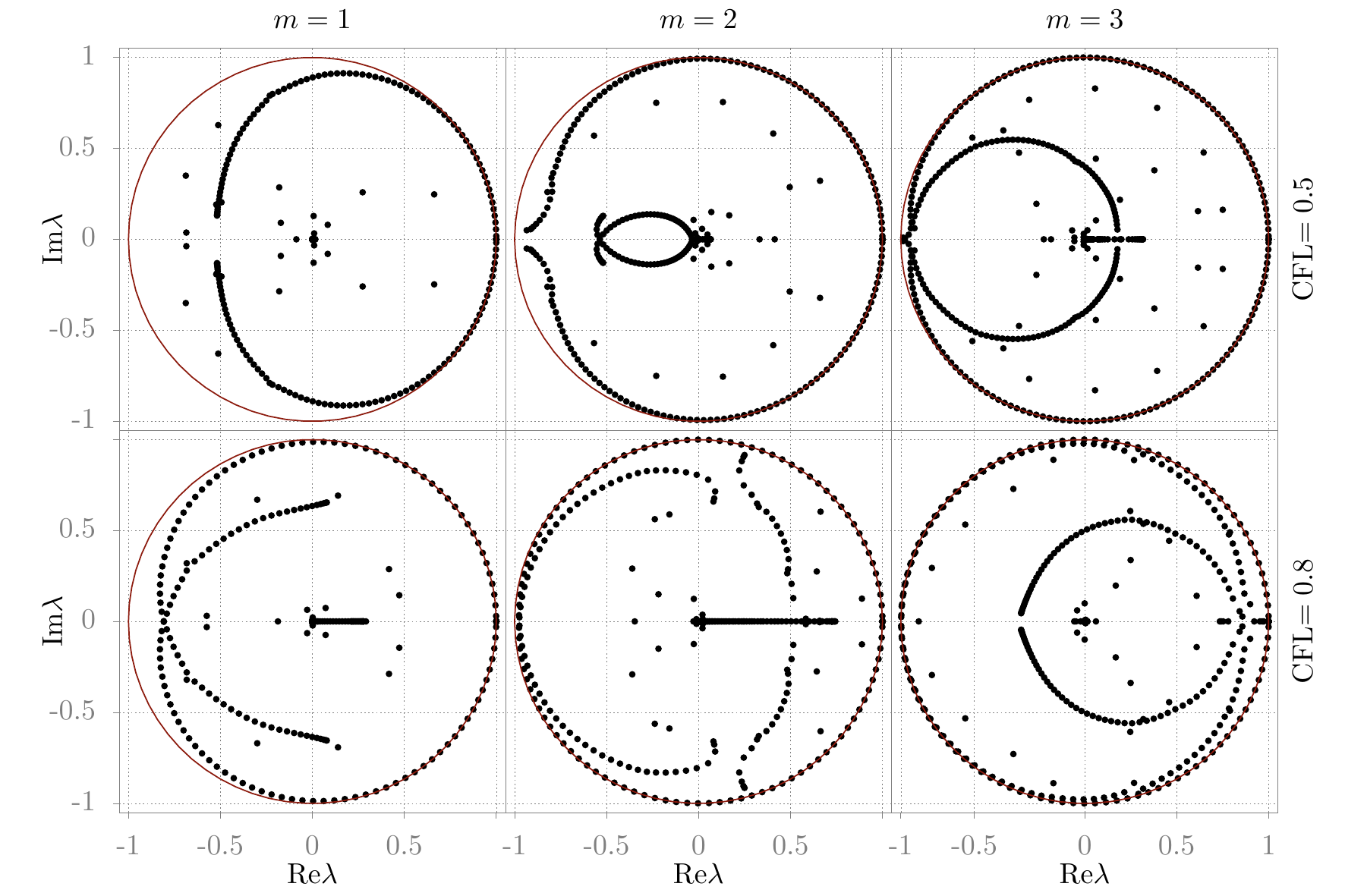}
    \caption{Spectrum of the amplification matrix $H$ for CFL numbers $\Delta t_a/h_a = 0.5, 0.8,$ orders of accuracy $3,5,7,$ and  \mbox{$n_a = 40,\ n_b = 5$}. No eigenvalues are outside the unit circle.\label{spectr2}}
    \end{center}
\end{figure}

Following the column-by-column construction process \eqref{constr} described \rr{above} we compute the amplification matrix $H$. The spectrum of $H$ is shown in Figure~\ref{spectr2} for $m=1,2,3$. Displayed results are for the cases  $n_a = 40$ and $n_b = 5$. The CFL numbers set for Hermite method are $\Delta t_a/ h_a = 0.5$ and $0.8$. The absolute value of eigenvalues do not exceed $1$. We note that if interpolation is used some eigenvalues of the amplification matrix shift outside of the unit circle.  Such unstable modes can possibly be stabilized by numerical dissipation / hyperviscosity but we do not pursue such stabilization here. Instead we observe that when projection is used all eigenvalues are inside the unit circle and the method is stable. Although we only display the results for one problem here the same results were obtained for other grid sizes, various overlap sizes to grid spacing ratios and different CFL numbers set for the Hermite method. We stress that it is possible to make the method unstable if we take the CFL number close to one and if we take $m$ to be larger than 3 and thus we only claim that the methods of orders of accuracy up to $7$ are stable.

\subsection{Convergence to an exact solution}
Using the same grid setup and boundary conditions as in the example above we test the method for the wave equation \eqref{eq1a}-\eqref{eq1b}, $c=1$ and initial conditions
\begin{align}
u(x,0) &= \sin\left(\frac{15\pi}2 x\right), \\
v(x,0) &= 0.
\end{align}
A solution to this problem is the standing wave
\be
u(x,t) = \sin\left(\frac{15\pi}2x\right)\cos\left(\frac{15\pi}2t\right). \nonumber
\ee
The errors for the solution on the grids are
\be
\varepsilon_a(x,t) = p_{i+{\hf}}(x,t) - u(x,t), x \in {x^a_i, x^a_i+1},\ i = 0, \dots, n_a,
\ee
for the Hermite grid and
\be
\varepsilon_b(x,t) = u^{h_b}(x,t) - u(x,t),
\ee
for the DG grid. The maximum error for the total method is
\be
\max \left( \max_{x \in {[x^a_0, x^{a}_{n_a}]}} |\varepsilon_a(x,t)|),\max_{x \in [x^b_0, x^b_{n_b}} |\varepsilon_b(x,t)|) \right). \nonumber
\ee

\begin{figure}[htb]
    \begin{center}
        \includegraphics[width=0.49\textwidth]{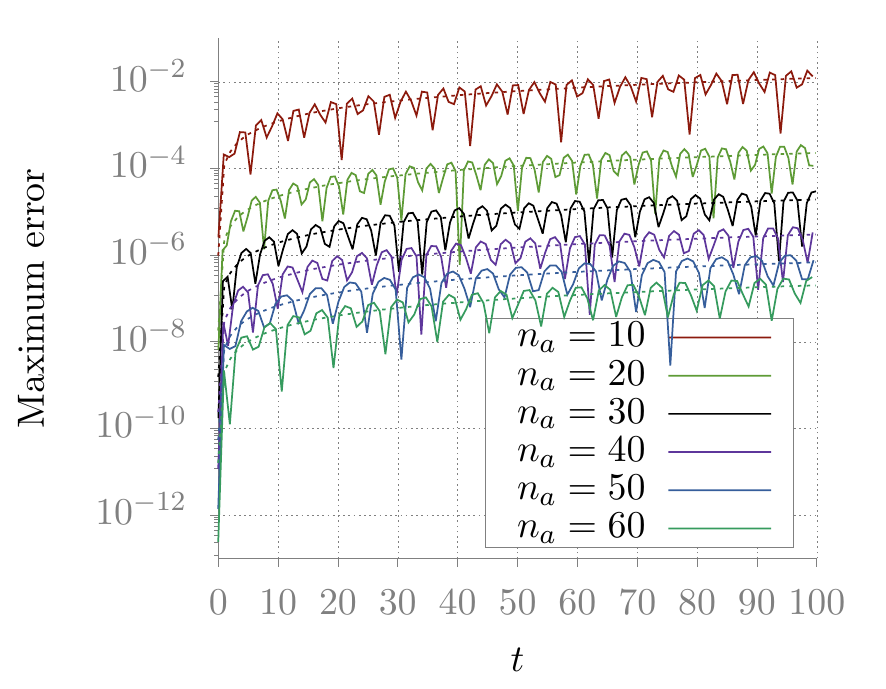}
        \includegraphics[width=0.49\textwidth]{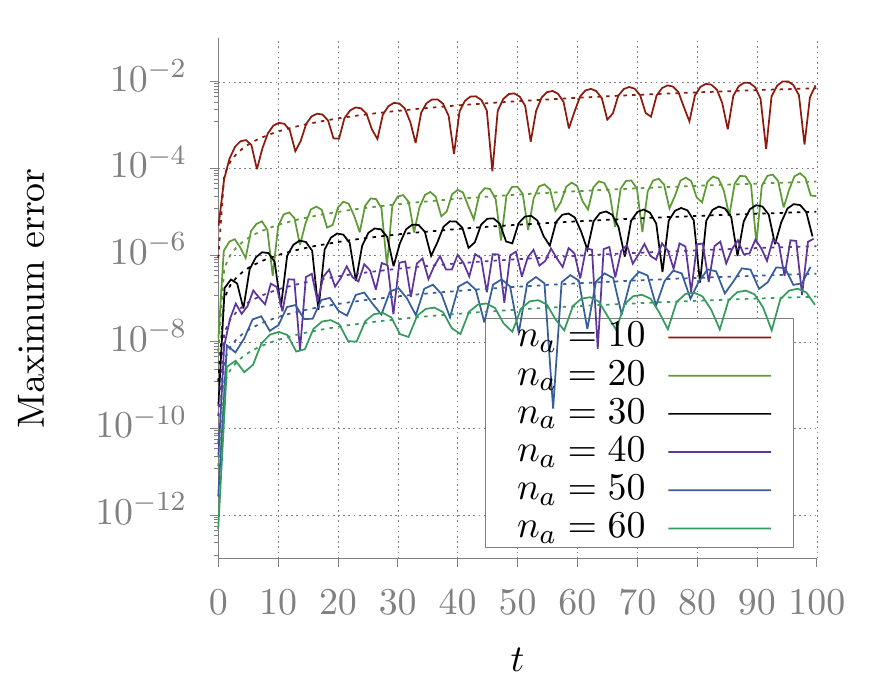}
        \caption{Maximum error of the solution as a function of time. The curves correspond to different refinements for $m$ = 3. In the left subfigure CFL number for Hermite method is set to $0.5$. In the right subfigure CFL number for Hermite method is set to  $0.75$. Dashed lines display lines $\alpha t$. \label{fe:1}}
\end{center}
\end{figure}

In Figure~\ref{fe:1} we display computed maximum errors as functions of time for the method with $m=3$ (i.e. the order of accuracy is 7). In the left subfigure the CFL number for the Hermite method is set to be $0.5$ and in the right subfigure the CFL number is set to be 0.75. For all Hermite grid sizes, the error growth is linear in time (dashed lines display a least squares fit of a linear function), indicating that the solution is stable for long time computations.

\begin{figure}[htb]
    \begin{center}
    \includegraphics[width=0.49\textwidth]{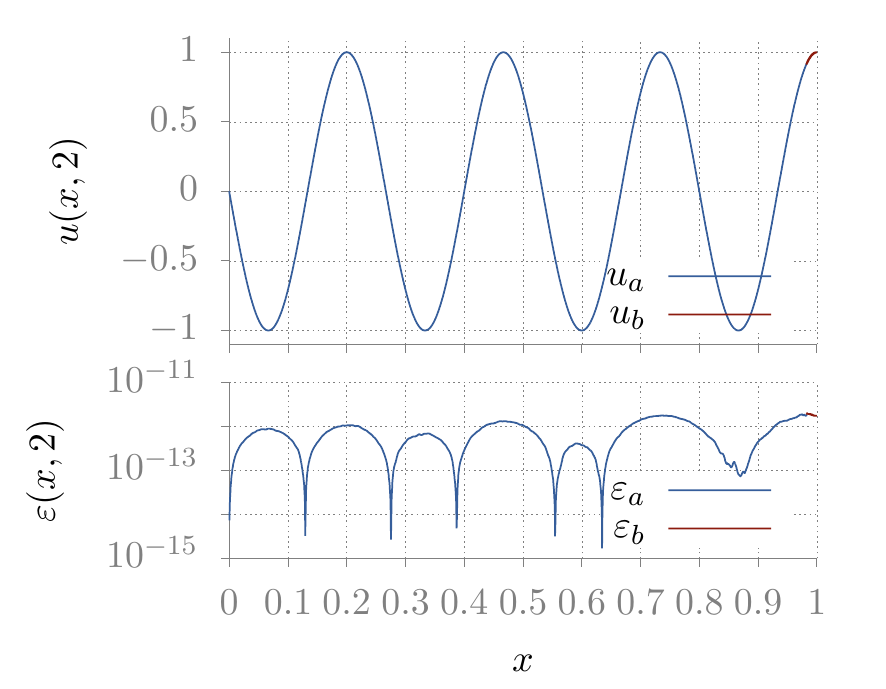}
    \includegraphics[width=0.47\textwidth]{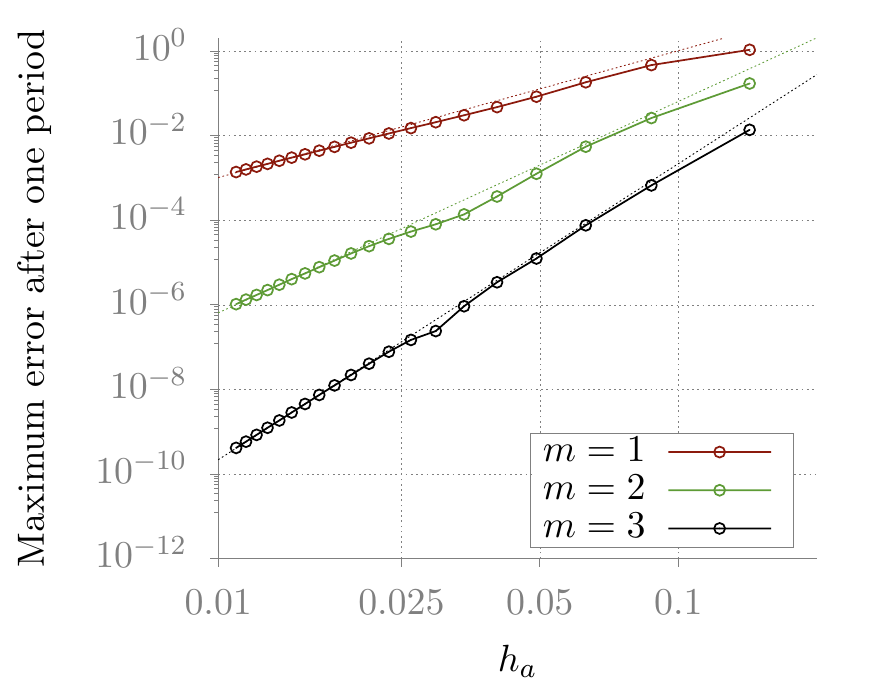}
    \caption{The upper left subfigure displays the solution at time $t=2$. The error at time $t=2$ is shown in the lower left subfigure. The number of grid points are $n_a = 200,\ n_b = 5$ and $m = 3$, the Hermite CFL number is set to $0.75$. Red curves indicate the solution and the error on the DG grid. Blue curves indicate the solution and the error on the Hermite grid. (The solution and the error were computed on finer grid, 10 grid points per cell/element). In the right subfigure we display a convergence plot for $m=1,2,3$. Dashed lines show the least squares fit of $C_mh_{a}^q,\ q = 3,5,7$.}
    \label{conv2}
\end{center}
\end{figure}

In the left subfigure of Figure~\ref{conv2} the numerical solution and the absolute error are shown for the $7$th order accurate method at time $t = 2$. As can be seen in the lower left subfigure in  Figure~\ref{conv2} the error is rather smooth across the overlap indicating that the projection is highly accurate.

To the right in Figure~\ref{conv2} the error at the final time $t=2$ is shown as a function $h = h_a$. The dashed lines show the least squares fit with polynomial functions of $h_a$ of order $3, 5$ and $7$ respectively. The results indicate that the orders of accuracy of the methods are $2m+1$ as expected. The parameters ($n_a$, $n_b$, $n_{DG}$, etc.) are the same is in previous example.

\subsection{Analytical solution in a disk. Rates of convergence}
\begin{figure}[hbt]
  \includegraphics[width=0.49\textwidth]{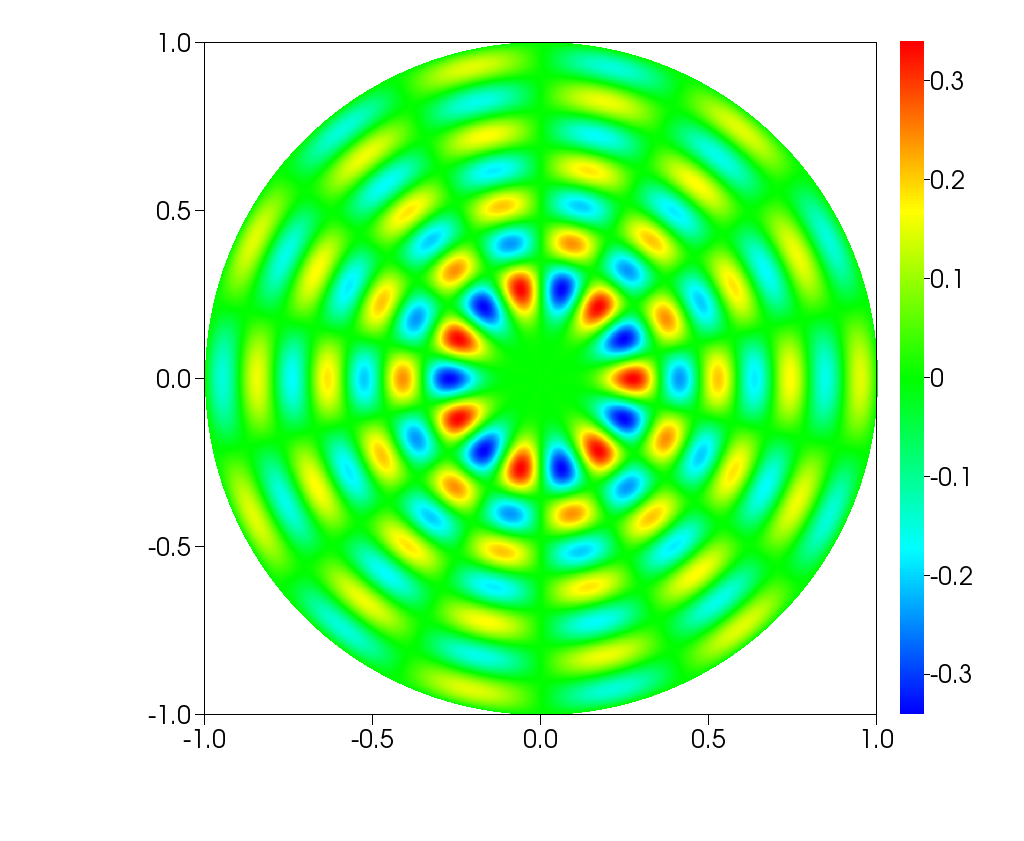}
  \includegraphics[width=0.5\textwidth]{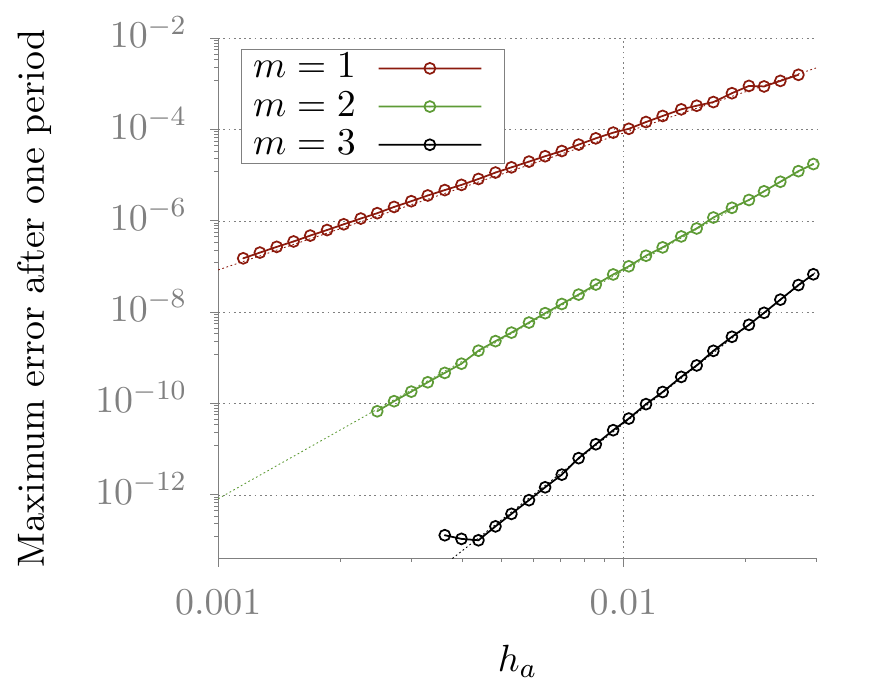}
  \caption{The left subfigure displays the initial condition.
  In the right subfigure the max-error at time $t = 2\pi/\kappa{77}$ as a function of grid spacing of the Hermite method. Solid curves correspond the methods with $m=1,2,3$ and dashed lines display the expected the convergence rates i.e. ${\cal O}(h_a^{2m+1})$.}
  \label{f:circ}
\end{figure}

Consider the solution of \eqref{eq1a}-\eqref{eq1b} with $f(x,y,t) \equiv 0$ on the unit disk, $(x,y) \in x^2+y^2 \le 1$, with homogeneous Dirichlet boundary conditions. Then the analytical solution can be expressed in polar coordinates as a composition of modes
\be
u_{\mu \nu}(r,\theta,t) = J_\mu(r\kappa_{\mu \nu})\cos(\mu \theta)\cos \kappa_{\mu \nu}t. \label{e:sol}
\ee
Here $J_\mu(z)$ is the Bessel function of the first kind of order $\mu$ and $\kappa_{\mu \nu}$ is the $\nu$th zero of $J_\mu$. In the following experiment we set $\mu=\nu=7$, $\kappa_{77} = 31.4227941922$. The initial condition $u_{77}(x,y,0)$ is displayed in the left subfigure of Figure~\ref{f:circ}.

\begin{figure}[htb]
  \includegraphics[width=0.49\columnwidth]{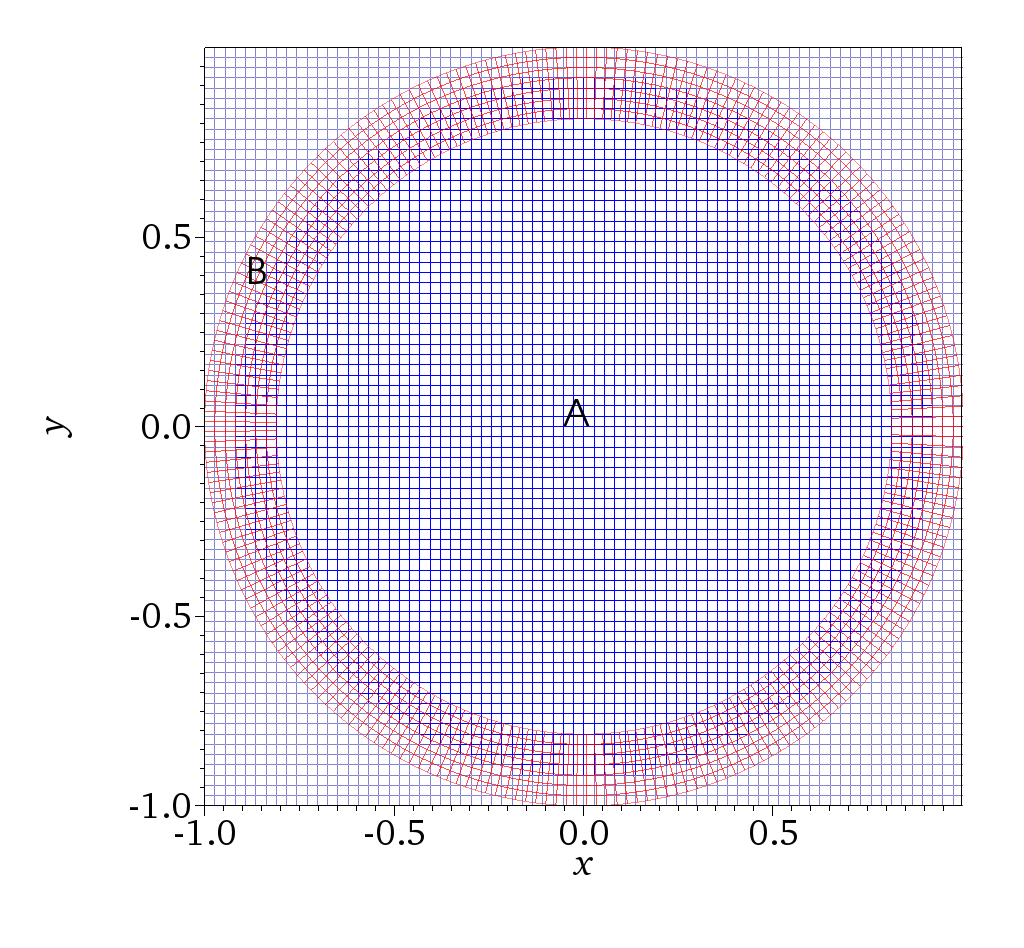}
  \includegraphics[width=0.49\columnwidth]{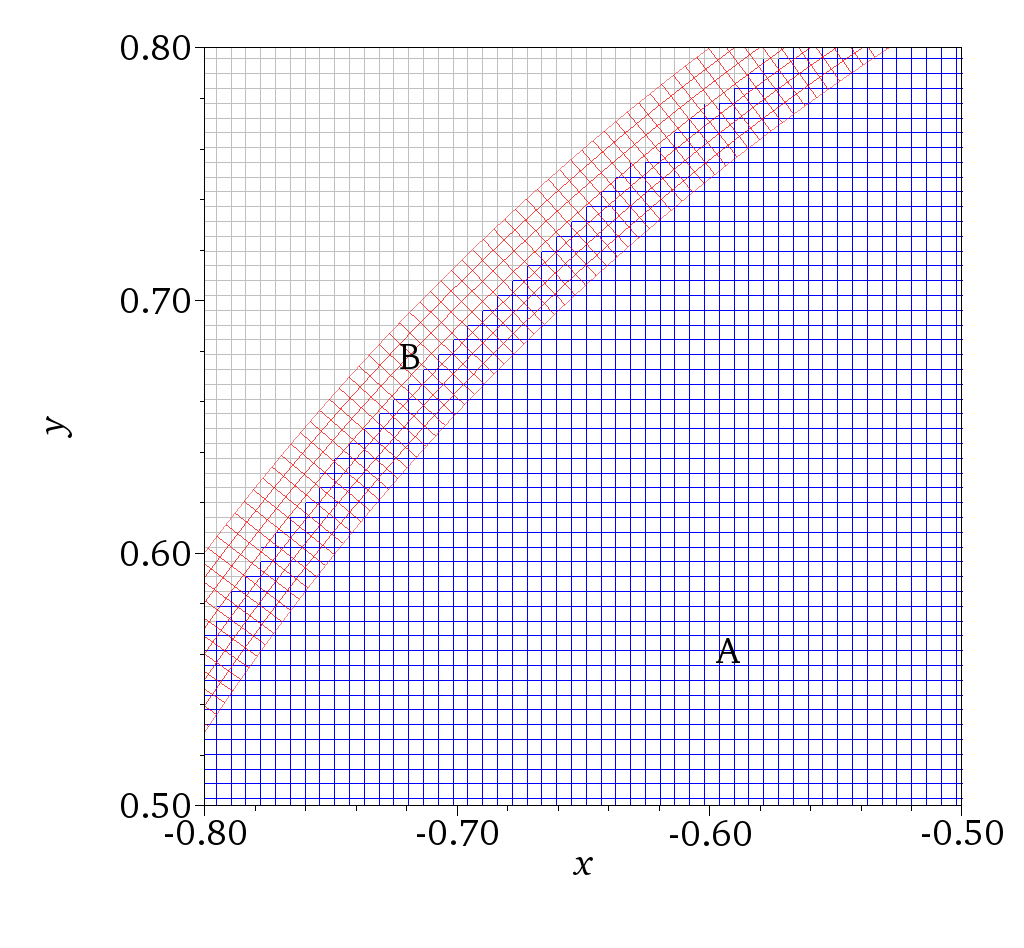}
  \caption{Overset grid set up for two different discretization widths. The Hermite grid is blue and DG grid is red. The Hermite grid is truncated at  radius $1-h_a$, i.e  one Hermite grid spacing smaller than the computational domain. This creates a stair sharped interior boundary. The solution at that boundary is imposed by by the projection described above. The curvilinear grid has 7 elements in the radial direction, thus the number of elements grows linearly with $n_a$. The number of grid points in the  Hermite grid grows as $n_a^2$.}
  \label{f:circ2}
\end{figure}

We setup overset grids as displayed in Figure \ref{f:circ2}. Grid $a$ is a Cartesian grid discretizing a square domain with $2n_a+1$ grid points in each direction and grid spacing $h_a = 1/n_a$. Grid $b$ is a curvilinear grid discretizing a thin annulus with radial grid spacing $1.1h_a$. For all refinements Grid $b$ has 7 elements in the radial direction thus the number of elements (or equivalently the number of DOFs of DG method) will grow linearly with the reciprocal of the discretization size $h_a$. In contrast the number of grid points in the Cartesian grid where the Hermite method will be used grows quadratically with $1/h_a$.

To measure the error we evaluate the solution on a finer grid, oversampled with 20 grid points inside each Hermite cell and DG element. The convergence is displayed in the right subfigure of Figure~\ref{f:circ}. The errors at time $t=2\pi/\kappa_{77}$ as  functions of $h_a$ for \mbox{$m = 1,2,3$} are displayed as solid lines. The dashed lines show the polynomials in $h_a$ of order $2m+1$. We use $h_a=1/34,1/36,...,1/94$ in the computations. As can be seen the expected orders of accuracy (3,5 and 7) are observed.
\begin{figure}[htb]
\begin{center}
    \includegraphics[width=0.39\textwidth]{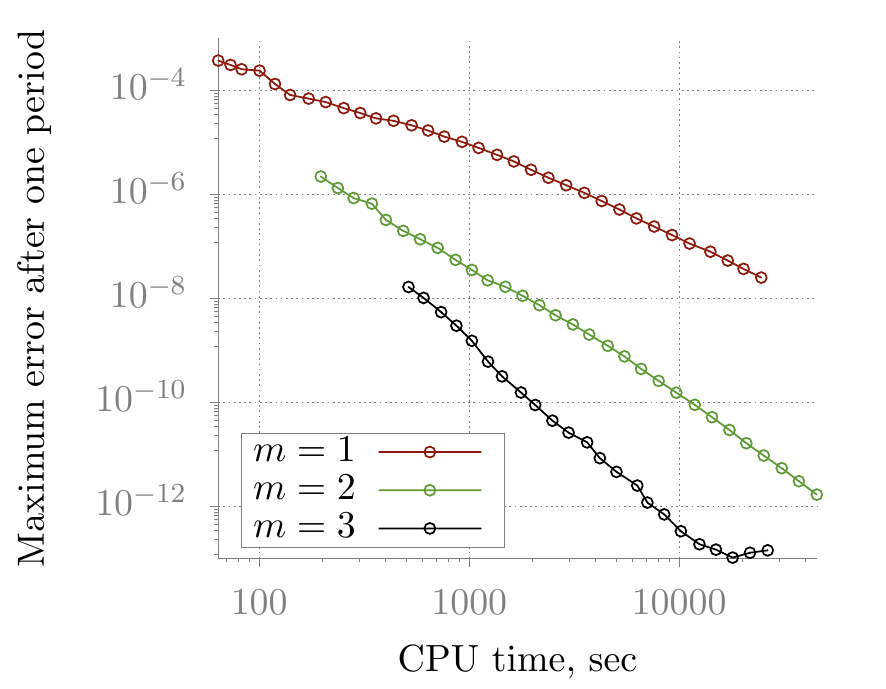}
    \includegraphics[width=0.39\textwidth]{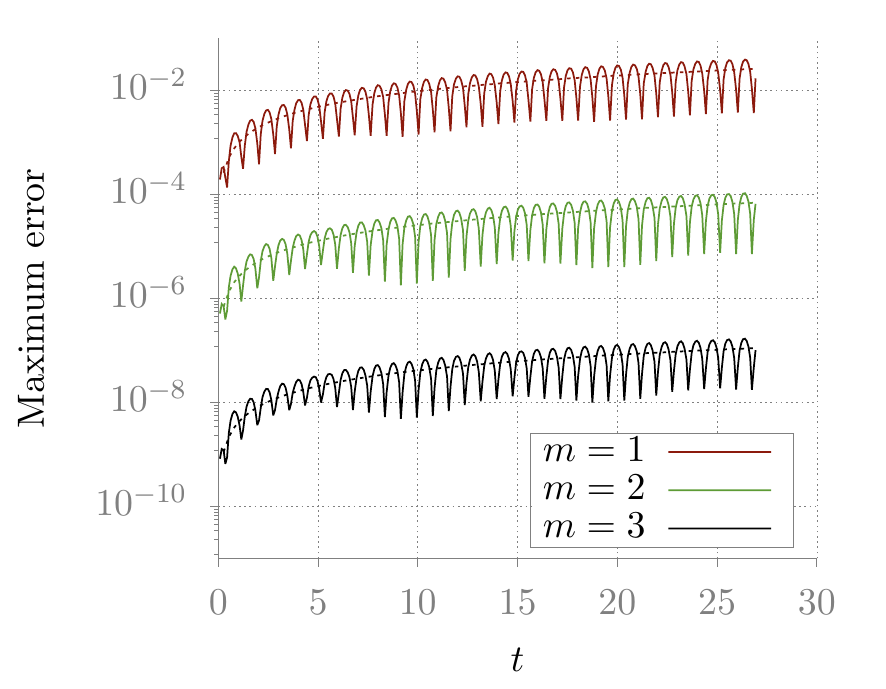}
    \caption{In the left subfigure the error at fixed time $t = 2\pi/\kappa_{77}$ is shown as a function of CPU time. The right subfigure displays the error as a function of time $t$.   Dashed lines are the linear functions $\alpha t$ formed by a least squares fit.}
  \label{f:circ3}
\end{center}
\end{figure}
To test the stability of the method we evolve the solution until time $t = 60\pi/7$ which is roughly $130$ periods of the solution. We set $h_a = 1/54$ and test methods with orders of accuracy $3, 5$ and $7$.
The error growth appears to be linear in time as indicated by dashed lines in the right subfigure of Figure~\ref{f:circ3}.

To test the performance of the method we evolve the method over one time period of the solution and measure the CPU time, see the left subfigure of Figure~\ref{f:circ3}. The red curve, displaying the error of the $3$rd, order accurate method only reaches the error $10^{-6}$ in about 1000 seconds while the $5$th and $7$th order accurate methods, using the same compute time, yield errors on the order of $10^{-8}$ and $10^{-10}$ respectively. Clearly the higher order methods are more efficient.

\rr{Table \ref{tab_speed1} displays a breakdown of time spent in the various parts of the code. As can be seen from the timing results the largest time is spent in the DG solver even for the finest grid. The increase in time does grow approximately quadratically and linearly for the Hermite and DG respectively so that eventually the complexity of the Hermite solver will dominate but practically speaking this may not happen for practical refinements for this problem. The large computational cost of the DG method is, in part, due to the small timestep requirement but also due to our implementation.}
\begin{table}[htb]
  \begin{center}
  \begin{tabular}{|c|c|c|c|c|c|}
    \hline
             &        HERMITE      &      DG  &      DG per step   &  H$\to$ DG    &     DG $\to$ H     \\
    \hline
    TIME     &            0.34     &    70.33 & 1.56289  &      10.00   &      11.32    \\
    \hline
     DOF     &          300448     &   222992 &   222992  &      15744   &      55748    \\
    \hline
     TIME / DOF   &    1.13(-6)         & 3.15(-4) & 7.00(-6)& 6.35(-4)     & 2.03(-4)      \\
    \hline
    \hline
TIME         &      1.61           &    159.67&    3.39  &   22.32     &    23.77      \\
    \hline
      DOF    &    1244760          &  451052  &  451052 &   32144      &     112763     \\
    \hline
      TIME / DOF   &  1.29(-6)           & 3.53(-4) & 7.51(-6) & 6.94(-4)     & 2.10(-4)    \\
    \hline
    \hline
      TIME      &      4.86    &    183.45  &  4.08 &    32.74    &     41.07    \\
    \hline
      DOF       &    5066944   &     905724   &   905724 &    64944   &     226431  \\
    \hline
      TIME/DOF    & 9.59(-7) &     2.02(-4)    & 4.45(-6)    & 5.04(-4) & 1.8(-4) \\
   \hline
  \end{tabular}
  \caption{Timing of the 7th order accurate hybrid Hermite-DG method for the disk experiment. The table contains timings for three different numbers of degrees of freedom. TIME denotes average time in seconds per 1 Hermite timestep of Hermite timestepping, DG timestepping and communication stages with the exception of the fourth column which displays the time per 1 DG timestep for the DG method. The TIME/DOF row in each block displays the time per degree of freedom computed by time evolution or communication. \label{tab_speed1}}
  \end{center}
\end{table}

\subsection{A wave scattering of a smooth pentagon}
In this experiment we study the scattering of a smooth pentagon in free-space. In addition to the use of non-reflecting boundary conditions experiment demonstrates the hybrid Hermite-DG method for the solution which is propagated over many wavelengths.
The geometry of the pentagon is defined as the smooth closed parametric curve:
\bea
x(s) &=& \frac 1{10}\left(1+\frac 1{10} \cos(10s)\right)\cos(s), \\
y(s) &=& \frac 1{10}\left(1+\frac 1{10} \cos(10s)\right)\sin(s), \ \ s \in [0,2\pi).
\eea
The pentagon  is placed in a square domain $(x,y) \in [-2,2]^2$ discretized by a Cartesian grid with grid spacing $1/n$, $n = 40$.
The curvilinear grid has 10 elements in the radial direction and the outer boundary is a circle of radius $0.1+20/n$. The overlap width is at most 5 DG elements.

On the boundary of the body we set Dirichlet data
\be
u(x,y,t) = \sin(\omega t),\ \ (x,y) \in \Gamma,\ \ t \ge 0,\ \ \omega = 250.
\ee
\begin{figure}
\begin{center}
  \includegraphics[width=0.39\textwidth]{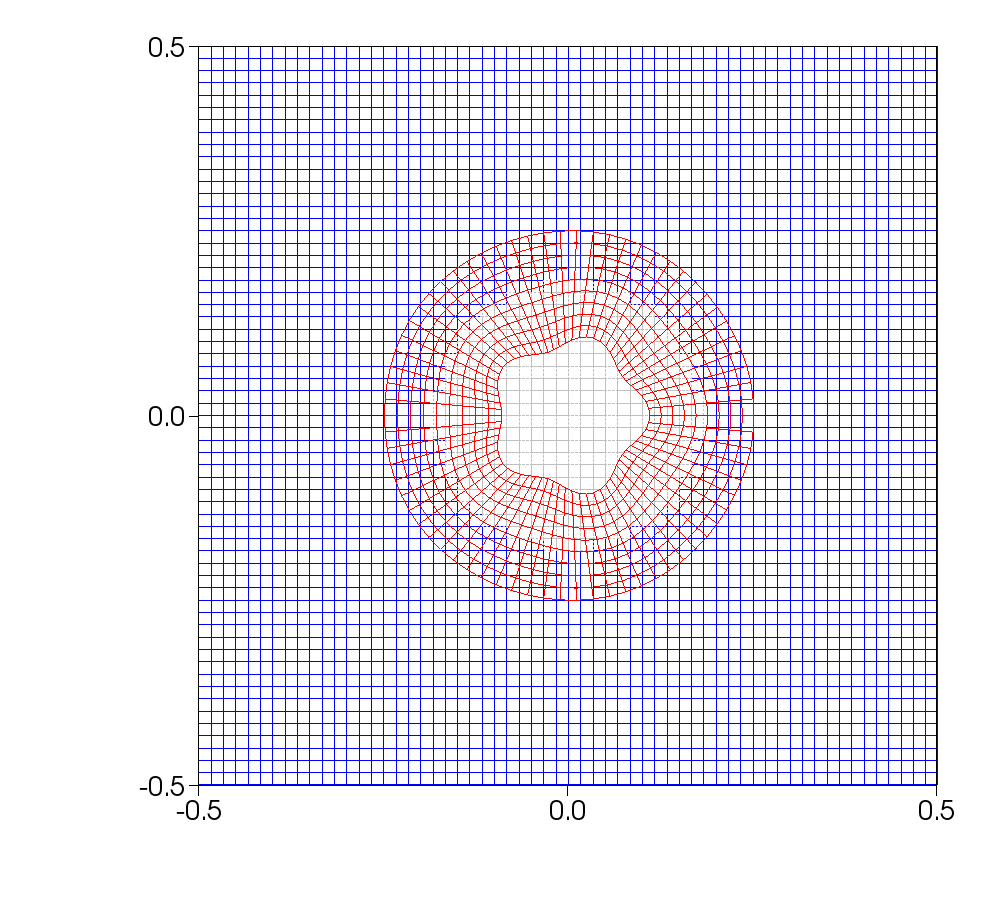}
  \includegraphics[width=0.4\textwidth]{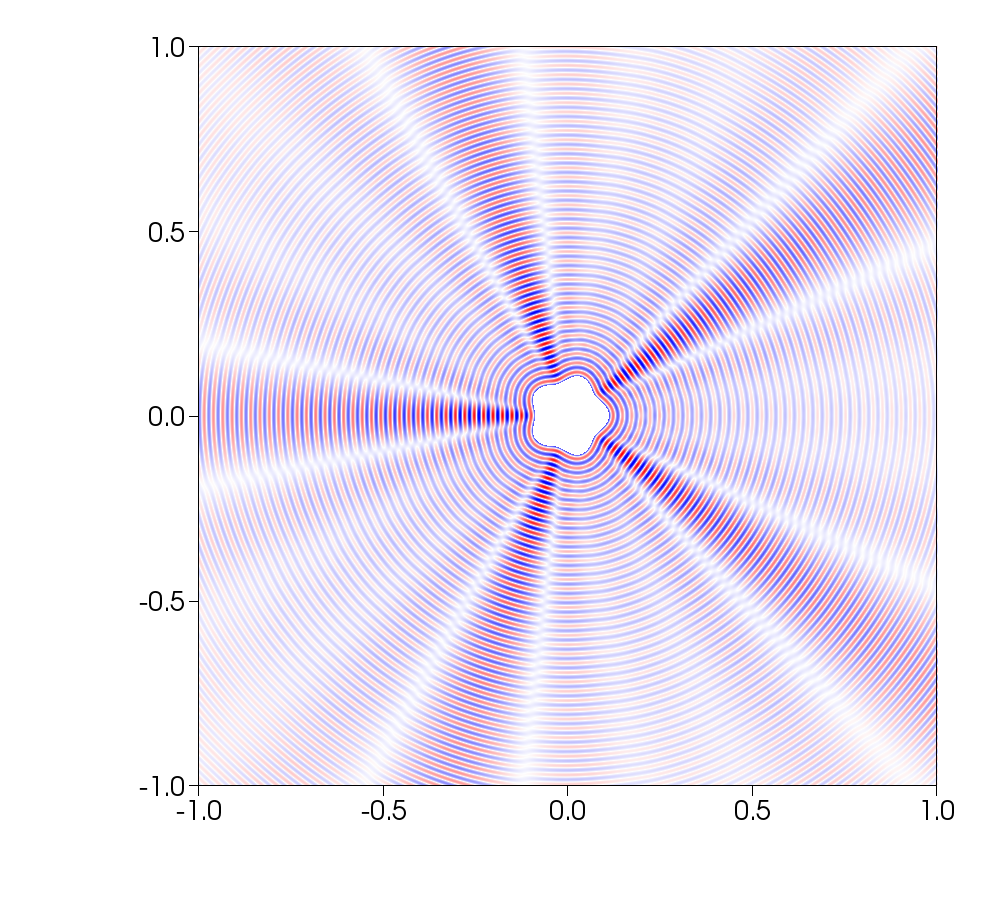}
  \caption{Left: Overset grid set up around the body. The overlapping DG grid and background Cartesian grid shown on the domain $[-0.5, 0.5]^2$.
  Right: Snapshot of $u(x,y,10)$.\label{fig:pent}}
  \end{center}
\end{figure}
The exterior boundary condition is modeled by truncating the domain using perfectly matched layers governed by the equations, (see \cite{appelo2012fourth} for derivation)
\be
u_{tt} = \frac{\partial}{\partial x}\left(
u_x + \sigma^{x} \phi^{(1)} \right)
+ \frac{\partial}{\partial y}\left(
u_y + \sigma^{y} \phi^{(2)} \right)
\sigma^{(x)}\phi^{(3)} + \sigma{(y)}\phi^{(4)}, \label{mod_pml}
\ee
where the auxiliary variables satisfy the equations
\be
\arraycolsep=1.4pt\def\arraystretch{2.2}
\begin{array}{r@{}l}
&\phi_t^{(1)}+(\alpha + \sigma{(x)})\phi^{(1)} = -u_x, \\
&\phi_t^{(2)}+(\alpha + \sigma{(y)})\phi^{(2)} = -u_y, \\
&\phi_t^{(3)}+(\alpha + \sigma{(x)})\phi^{(3)} = -u_{xx}-\frac{\partial}{\partial x}
\left(\sigma^{(x)}\phi^{(1)} \right), \\
&\phi_t^{(4)}+(\alpha + \sigma{(y)})\phi^{(4)} = -u_{yy}-\frac{\partial}{\partial y}
\left(\sigma^{(y)}\phi^{(2)} \right).
  \label{newvars}
\end{array}
\ee
The damping profiles $\sigma^{(z)}$, $z = x,y$ are taken as
\be
\sigma^{(z)}(z) = \sigma_s \left(\tanh\left(\frac{z-z_1}{0.7w_{\rm lay}}\right)
- \tanh \left(\frac{z-z_2}{0.7w_{\rm lay}}\right)\right). \nonumber
\ee
\begin{figure}
\begin{center}
  \includegraphics[width=0.4\textwidth]{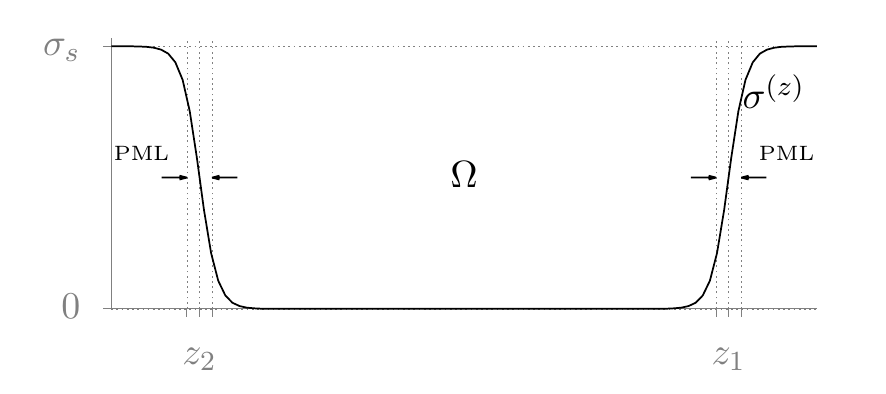}
  \caption{The damping profile $\sigma^{(z)}$ for PML with damping strength $\sigma_s$. The damping is zero at the center of the domain $\Omega$ and rapidly increases in the PML on the both edges of the domain.}  \label{sig}
\end{center}
\end{figure}

Here $\sigma_s$ is a damping strength, $w_{\rm lay}$ is layer width and $z_1$ and $z_2$ control the location of the damping profile of the PML. The shape of $\sigma^{(z)}$ is displayed in Figure~\ref{sig}. In experiments involving PML we discretize the modified equations \eqref{mod_pml} with the Hermite method.
The order of accuracy of the methods is set to be 7, i.e. $q_u = q_v = 6$, and $m = 3$.  The solution is evolved to $t=10$. A snapshot of the solution at the final time is displayed in the right subfigure of Figure~\ref{fig:pent}. The proposed algorithm clearly is able to accurately propagate waves in complex domains.

\begin{table}[]
  \begin{center}
  \begin{tabular}{|c|c|c|c|c|c|}
    \hline
             &        HERMITE      &      DG  &      DG per step   &  H$\to$ DG    &     DG $\to$ H     \\
    \hline
    TIME &  23.87 & 14.15  & 0.25 & 1.98 & 0.27 \\
    \hline
     DOF &       1972100 &        87010  &       87010 &         3280 &         7910 \\
    \hline
     T/DOF &  0.12(-4) & 0.16(-3)  &0.29(-5) & 0.60(-3) & 0.34(-4) \\
    \hline
    \hline
   TIME  &  40.10 & 17.52  & 0.37  & 2.52  & 0.42 \\
    \hline
     DOF &       4170520 &       125543  &      125543  &        4920  &       11413 \\
    \hline
     T/DOF &  0.96(-5) & 0.13(-3)  &0.30(-5)  &0.51(-3)  &0.37(-4) \\
    \hline
    \hline
    TIME &  76.18 & 27.38  & 0.60  & 3.94  & 0.70 \\
    \hline
     DOF &      11007352 &       203852  &      203852  &        8200  &       18532 \\
    \hline
     T/DOF &  0.69(-5) & 0.13(-3)  & 0.29(-5)  & 0.48(-3)  &0.38(-4) \\
    \hline
    \hline
    TIME &  162.28 &  43.08  & 0.87  & 7.30  & 1.38 \\
    \hline
     DOF &      34433112 &       287811  &      287811  &       15088  &       31979 \\
    \hline
     T/DOF &  0.47(-5) & 0.14(-3)  &0.31(-5)  &0.48(-3)  &0.43(-4) \\
    \hline
  \end{tabular}
  \caption{Timing of the 7th order accurate hybrid Hermite-DG method for the smooth pentagon experiment. The table contains timings for three different numbers of degrees of freedom. TIME denotes average time in seconds per 1 Hermite timestep of Hermite timestepping, DG timestepping and communication stages with the exception of the fourth column which displays the time per 1 DG timestep for the DG method. The T/DOF row in each block displays the time per degree of freedom computed by time evolution or communication.}
  \label{tab_speed2}
  \end{center}
\end{table}

\rr{Table \ref{tab_speed2} displays a breakdown of time spent in the various parts of the code. As can be seen from the timing for this problem the largest time is now spent in the Hermite solver. Here, due to the geometry being an interior object, the relative number of degrees of freedom in the DG solver is small and we see the asymptotic behavior more clearly than for the disc experiment. }

\subsection{Wave scattering of many cylinders in free space}
As another demonstration of the method we simulate a domain with multiple circular holes. Precisely we consider the infinite domain $\Omega \in [-\infty,\infty]\times[-\infty,1.33]$ with homogeneous Neumann boundary condition at $y = 1.33$. The computational domain is a rectangle $[-1,1]\times[-1.33]$ with PML $|x| > 1$ and $y < -1$. Inside the computational domain there are 5 cylinders of radii $0.1$ and centers at $(x_k,y_k), k = 1,\dots 5$. We impose the homogeneous Dirichlet boundary conditions on the boundary of all cylinders except first. On the first cylinder we impose a time dependent boundary condition
\be
u(t,x,y) =
(t-0.1)\exp\left(-918(t-0.1)^2\right),\ (x,y) \in \{ (x-x_1)^2 + (y - y_1)^2 = 0.01\}.
\nonumber
\ee
The initial solution is at rest.

\begin{figure}[]
  \begin{center}
    \includegraphics[width=.245\textwidth]{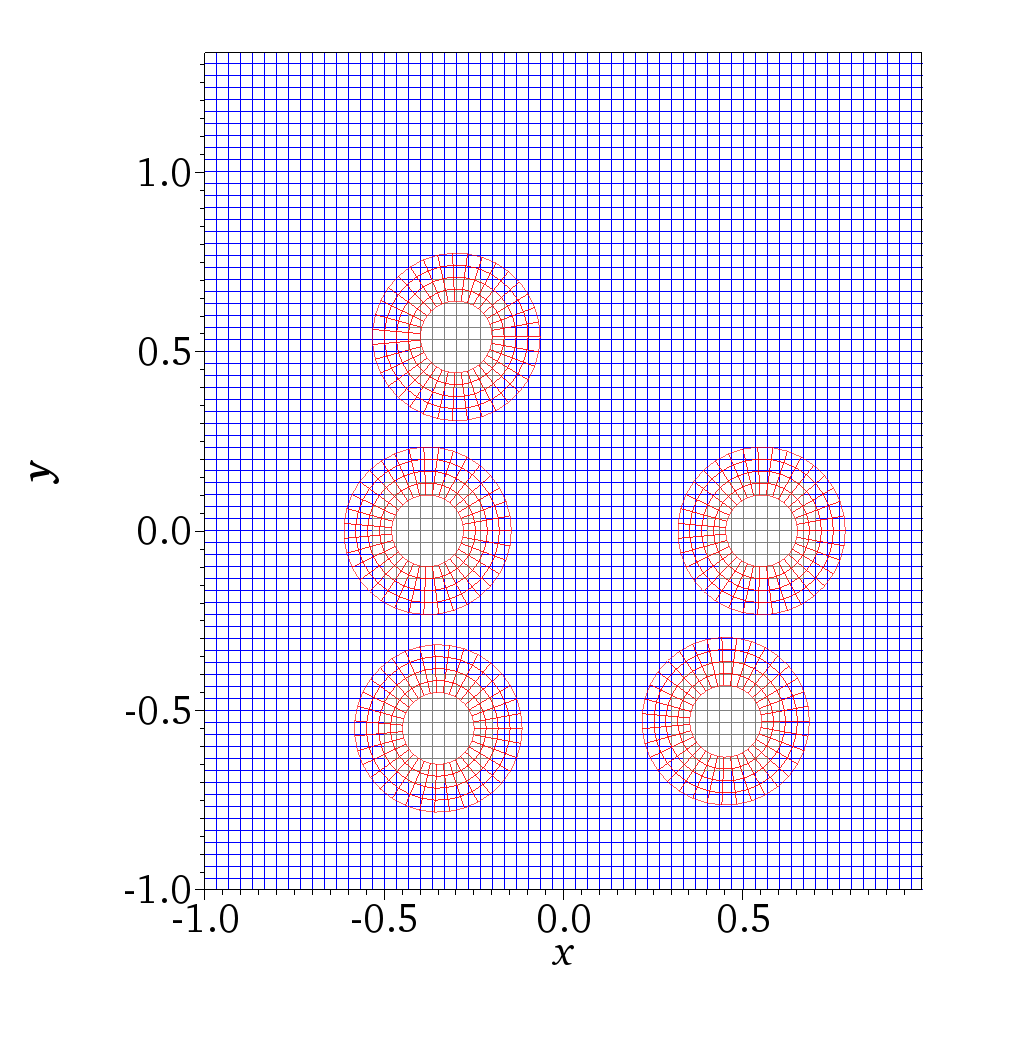}
    \includegraphics[width=.245\textwidth]{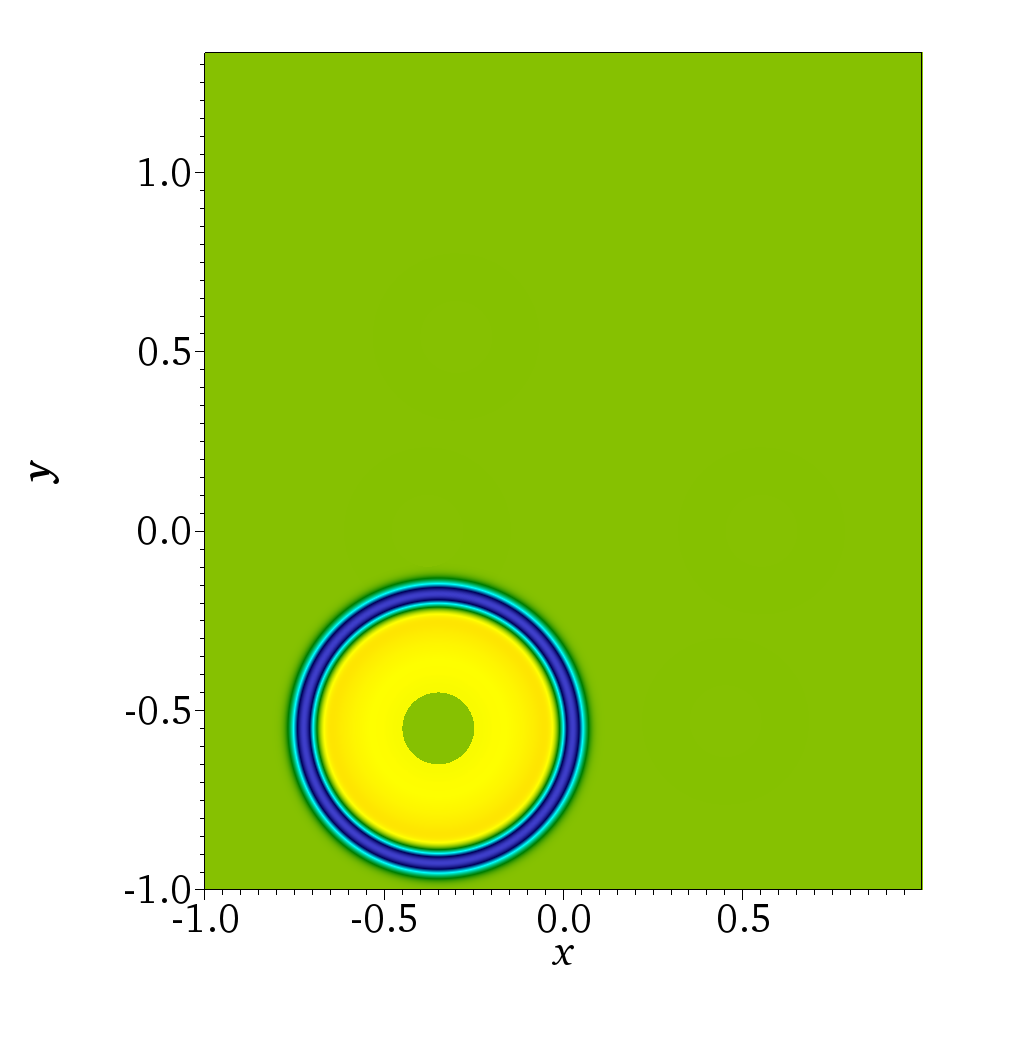}
    \includegraphics[width=.245\textwidth]{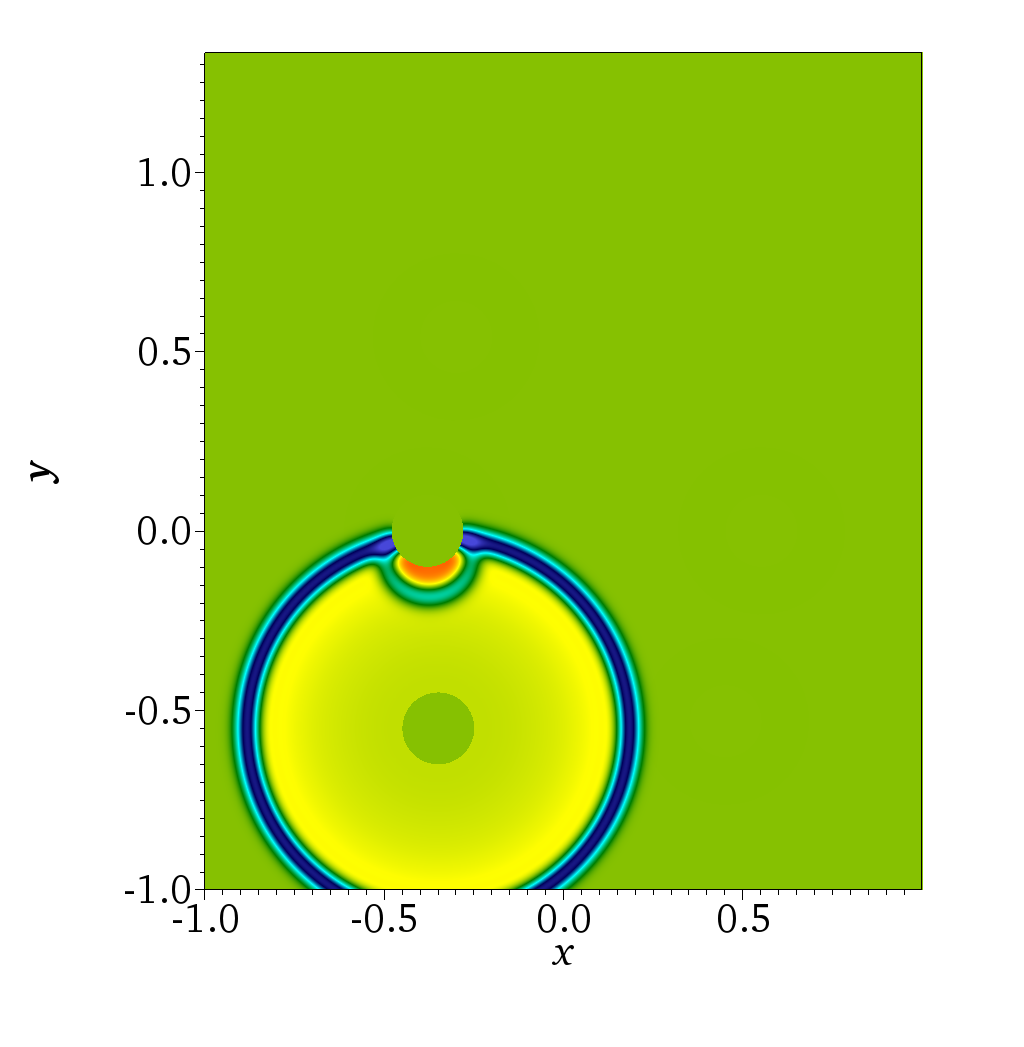}
    \includegraphics[width=.245\textwidth]{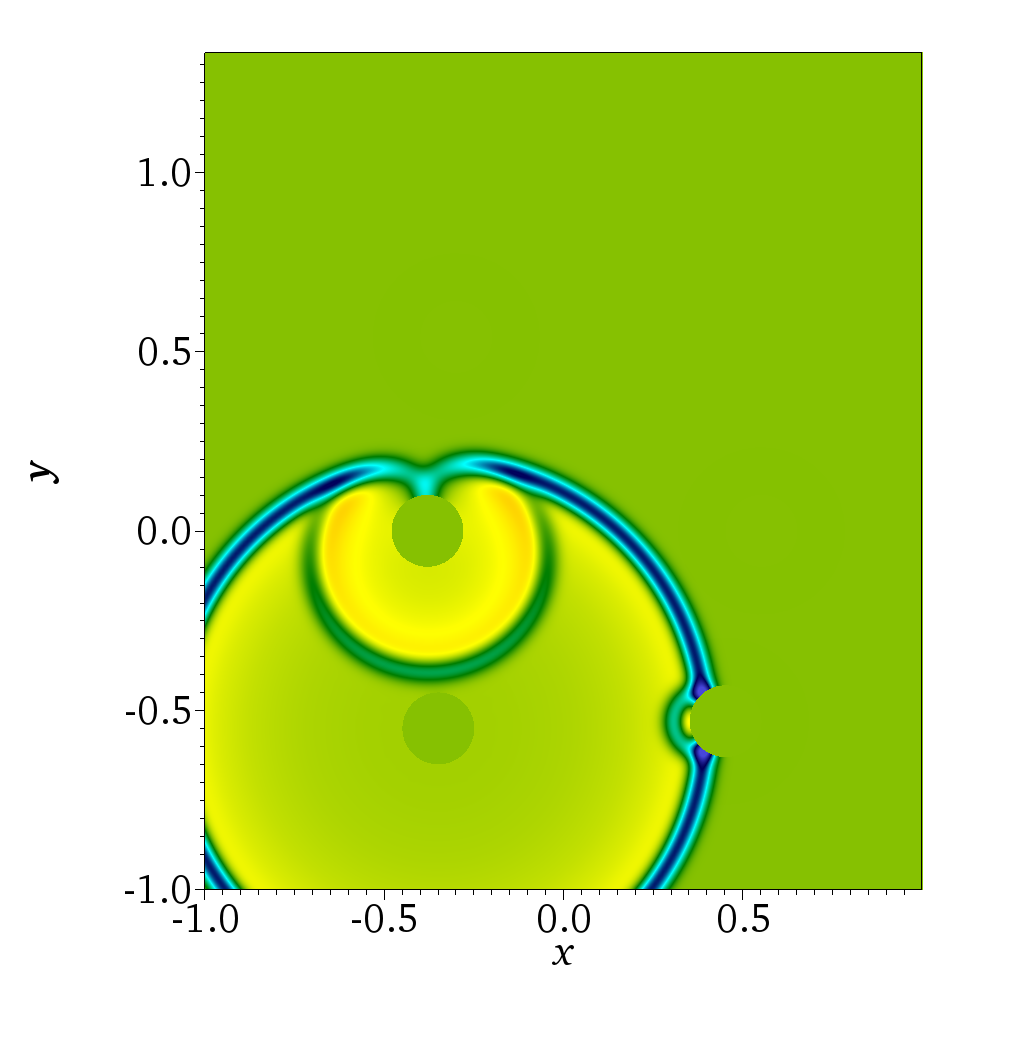}
    \includegraphics[width=.245\textwidth]{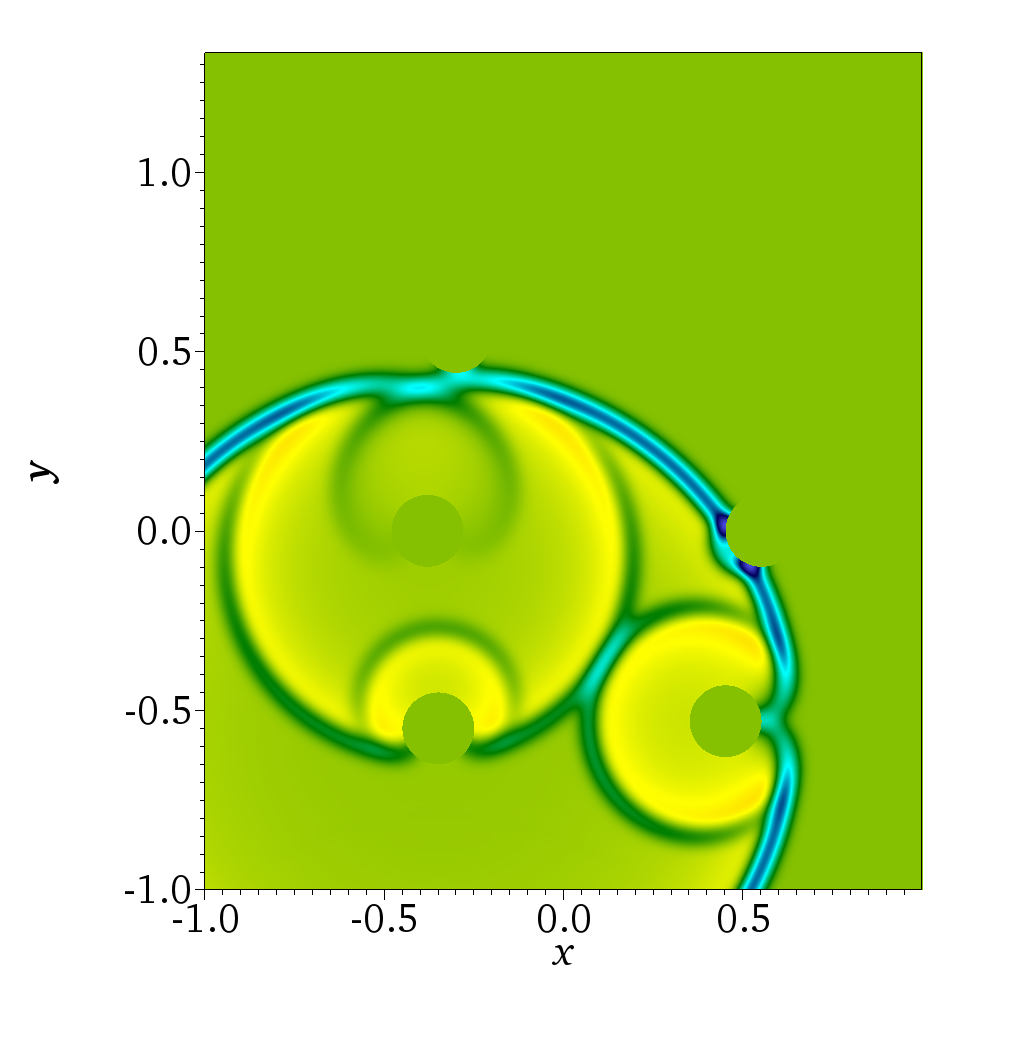}
    \includegraphics[width=.245\textwidth]{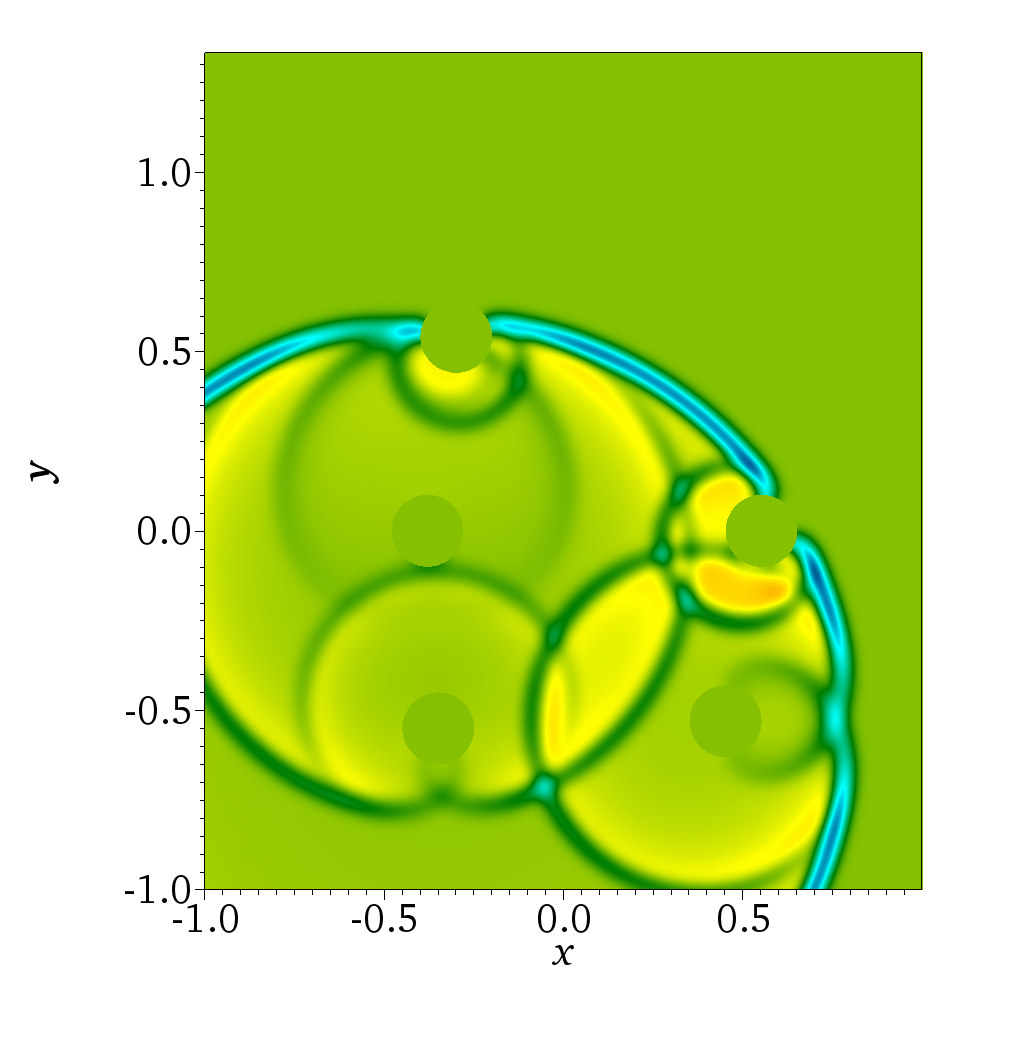}
    \includegraphics[width=.245\textwidth]{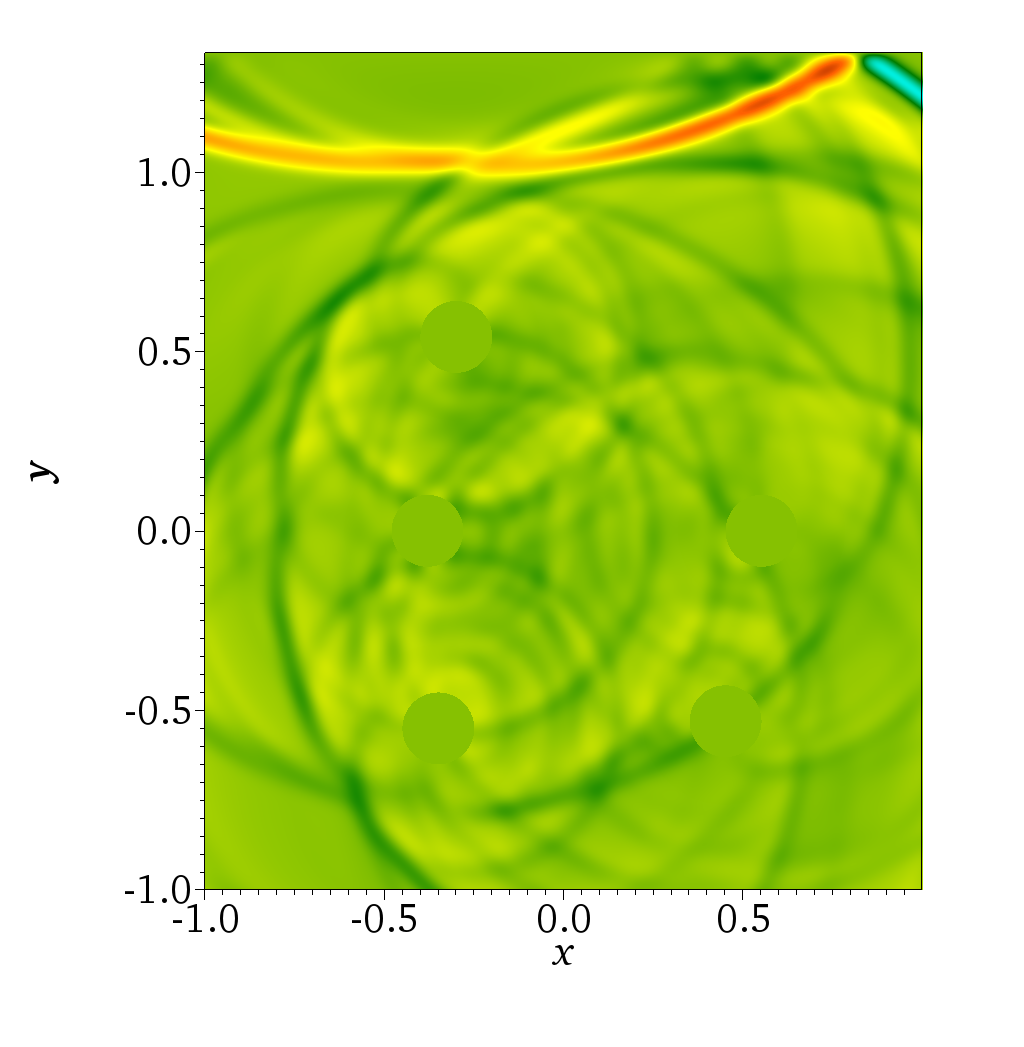}
    \includegraphics[width=.245\textwidth]{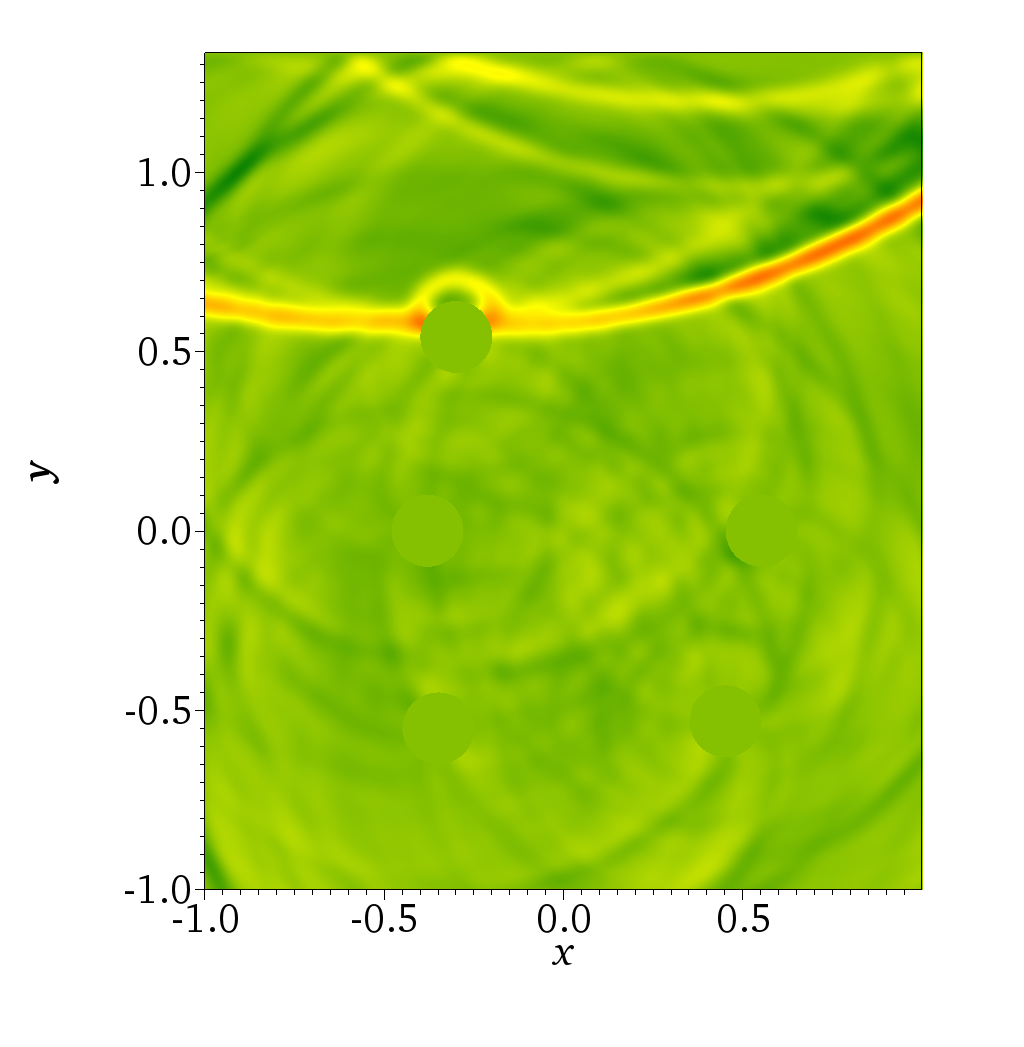}
\caption{Overset grid setup and solution plots for 5 bodies in a free half space. In the upper left subfigure the grids are shown: DG gridlines plotted with red color, Cartesian grid lines inside the domain plotted with blue. Other figures are the solution plots at various increasing times.}
  \label{fig:5b}
  \end{center}
\end{figure}

The set up of the numerical method is similar to the  previous experiment. The Cartesian grid covers the background domain and the PML. The 5 circular grids are placed around the bodies as shown in the upper left subfigure Figure~\ref{fig:5b}. In this experiment we used the 7th order method. It can be noticed that as in all solution plots provided in this paper the solution is smooth across the overlap due to the high accuracy of methods used and the projection used for communication.

\subsection{An inverse problem, locating a body in free space} \label{ss:body}
As a final experiment we solve the inverse problem of locating a cylindrical body in free space. An application of this problem could be a to locate a tunnel under the ground and determining its radius by sending waves from source devices buried at a relatively small distance from the surface and recording the solution near the surface. Waves will propagate from a source, reflect from an underground cavity and travel back to the surface to be captured by the recording devices. The underground cavity can be located by minimizing a cost functional, i.e misfit function of recorded data and data obtained from the numerical simulation in each iteration of the optimization process.

\begin{figure}[htb]
  \begin{center}
\includegraphics[width=.39\textwidth]{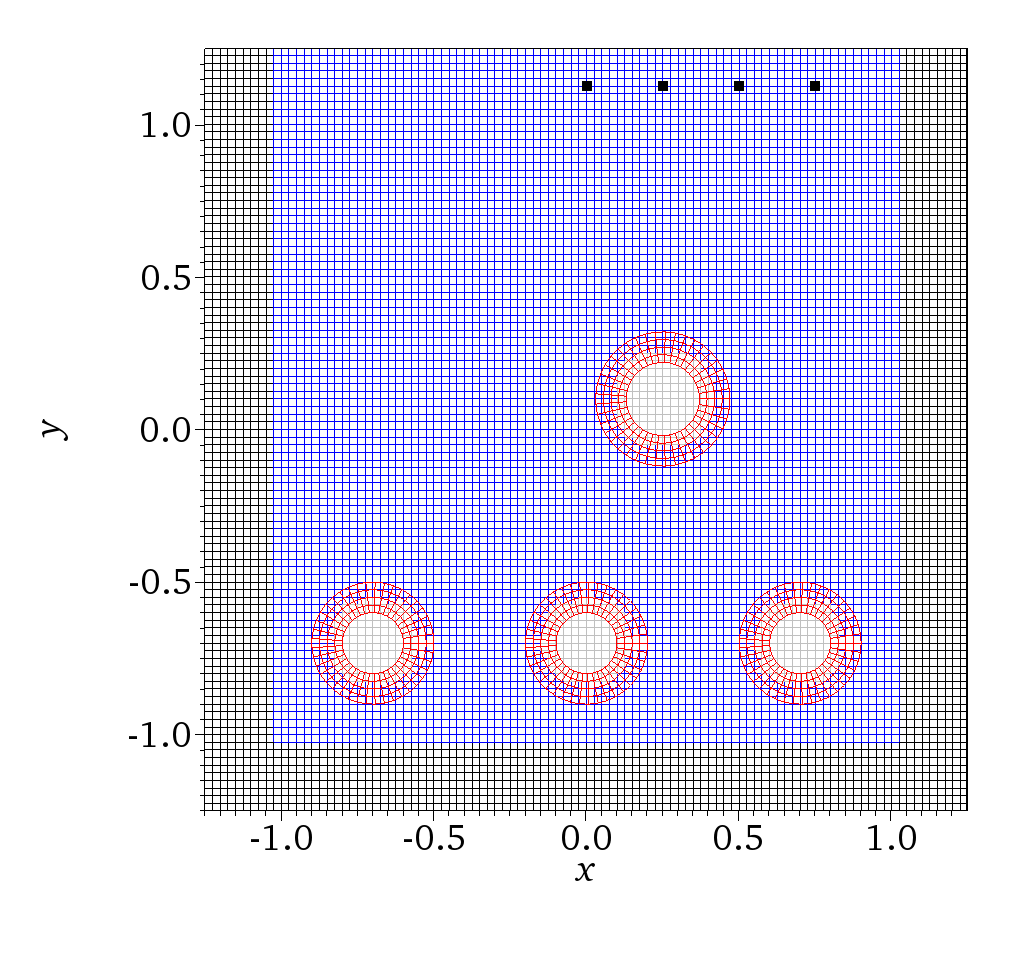}
\includegraphics[width=.39\textwidth]{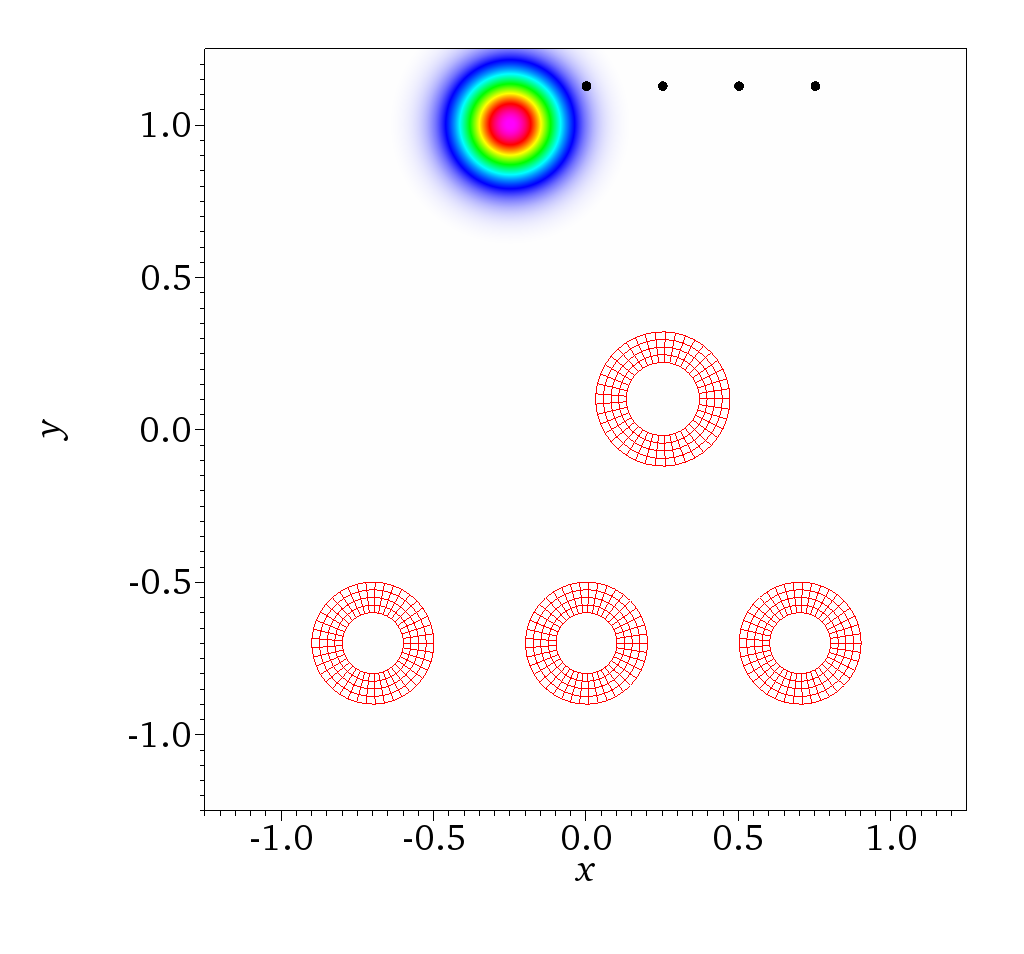}
\caption{The inverse problem set up. Receivers are marked as dots. The left subfigure displays the complete overset grid set up, with 4 DG grids around bodies,
and a Cartesian background grid. Blue and red grids discretize physical subdomain; the black grid is the PML layer; the gray grids indicate the subdomains covered by DG grids.
The right subfigure displays of the initial condition, a smooth Gaussian centered at (-0.25, 1), receivers and the DG grids.  \label{fig:inv}}
  \end{center}
\end{figure}

Consider a square region $\Omega \in [-1,1]\times[-1,1.25]$ with 3 circular bodies of radius $r = 0.1$ with centers at $x_1 = -0.7$, $x_2 = 0$, $x_3 = 0.7$ and $y_1=y_2=y_3 = -0.7$. On the boundary of the bodies we impose homogeneous Dirichlet boundary conditions. On the top boundary $y = 1.33$, that acts as a "ground surface" we impose homogeneous Neumann boundary conditions. The exterior boundary conditions at $x=\pm 1$ and $y = -1$ are imposed by truncating the domain using PML. We discretize the domain with a Cartesian grid. Around each of the cavities we place annular DG grids that are 5 cell wide. An example of a complete set up with $4$ receivers is shown in Figure~\ref{fig:inv}.

First we create synthetic data by recording the displacement $u$ at equidistant locations of the receivers \[(0,0.125),\ (0.25,0.125),\ (0.25,0.125),\ (0.25,0.125),\] to time $T=2$. Let there be another circular body of radius $A_1$ and center at $(x,y) = (A_2,A_3)$ we want to locate. In the right figure the first source is active, i.e. the initial condition is a smooth Gaussian centered at $\hat x = -0.25$, $\hat y = 1$, with no initial velocity
\be
u(0,x,y) = \exp\left(-40\left((x-\hat x)^2+(y-\hat y)^2\right)\right),\ \  v(0,x,y)\equiv 0. \label{in:gauss}
\ee

First we create the synthetic data for the exact location of the target, $A_1^\ast = 0.12$, $A_2^\ast = 0.25$ and $A_3^\ast = 0.1$. This gives us $u^\ast$. To locate the cavity we minimize the cost function that is a sum of squared $L_2$ norms of discrepancies between the output of the numerical simulation and synthetic data $u^*(t,\check x_l, 0.125)$
\be
F(A_1, A_2, A_3) = \sum_{l = 1}^4 \int_0^T \left(u(t,\check x_l, 0.125) - u^*(t,\check x_l, 0.125) \right)^2.
\ee

During the minimization we impose the bounds $0.01 \le A_1 \le 0.2$, $|A_2| < 0.5$, $|A_3| < 0.2$. To recover $A^\ast_1, A^\ast_2$ \bb{and $A_3^\ast$} we use the L--BFGS-B algorithm, (see \cite{byrd1995limited} for description). Forward differences are used to compute the gradients, resulting in total $1+3$ simulations per iteration. Table~\ref{tab1} displays the convergence results in detail for the initial values \bb{at $1\%$ of the exact solution,  that are $0.101$, $0.2525$ and $0.1212$ respectively}. At the $6$th iteration the values computed were $0.1$, $0.25$ and $0.12$, accurate to the $8$th digit. \bb{For the initial guesses with larger the $1\%$ deviation from the exact solution, it becomes harder to converge to a global minimum. The minimization process would become more robust if more data is recorded at the receivers, for example by increasing the number of receivers, recording longer data traces or adding simulations with different initial conditions}.

\bb{Although during the minimization process before each simulation the grids have to be regenerated this is inexpensive since the grid generation is local. Precisely in each new iteration the DG grid is adjusted to by regenerating an annular grid based on the updated center location and radius.
}

\begin{table}
  \begin{center}
  \begin{tabular}{|c|c|c|c|c|c|}
    \hline
    $N$ iter. & $A_1$ & $A_2$ & $A_3$ & $F$ & $\|\nabla F\|$ \\
    \hline
    $0$ & $0.10100$ &  $0.25250E$ &  $0.12120$ &$7.83582(-6)$ & $4.41312(-3) $ \\
    \hline
    $1$ & $0.10117$ &  $0.25077E$ &  $0.11944$ &$1.96547(-7)$ & $2.88478(-4) $ \\
    \hline
    $2$ & $0.10117$ &  $0.25066E$ &  $0.11955$ &$1.38952(-7)$ & $2.32250(-4) $ \\
    \hline
    $3$ & $0.10105$ &  $0.25012E$ &  $0.12000$ &$2.26171(-8)$ & $4.32507(-5) $ \\
    \hline
    $4$ & $0.10095$ &  $0.25010E$ &  $0.12001$ &$1.87633(-8)$ & $3.97138(-5) $ \\
    \hline
    $5$ & $0.10002$ &  $0.24999E$ &  $0.12001$ &$5.81672(-11)$ & $5.80355(-6) $ \\
    \hline
    $6$ & $0.10000$ &  $0.25000E$ &  $0.12000$ &$1.00121(-15)$ & $4.21271(-8) $ \\
    \hline
  \end{tabular}
  \caption{Convergence results of L-FBGS-B algorithm for the inverse problem for locating a body in free space. At each iteration the cost function $F$ and its gradient $\nabla F$ is computed from the numerical solution of the wave equation. The forward solver is implemented using the $5$th order accurate Hybrid Hermite-DG overset grid method.}
  \label{tab1}
  \end{center}
\end{table}

\section{Summary}
We have presented overset high order numerical methods for numerical solution of the wave equation. The hybrid H--DG overset grid method combines the highly efficient Hermite method  on Cartesian grids with a DG method to treat  complex boundaries. To combine the methods the overset grids were used. The advantage of using the overset grids for complex boundary problems is the low computational cost that asymptotically approaches the cost of the  Cartesian solver.

In this work we communicate solutions via $L_2$ projection and this procedure combined with the dissipative nature of the methods was observed to be sufficient to guarantee stability without the need to add any artificial dissipation.

Stability, accuracy and efficiency of the method were tested numerically. To test the stability in 1 dimension, we looked at the spectrum of the amplification matrix associated with the method. For CFL numbers $< 0.75$ for the Hermite method, the overall method was stable in all tested settings for grid sizes and orders of accuracy $3, 5$ and $7$. In 1 and 2 dimensions we also tested the stability by displaying the error growth as a function of time for long times.

Finally,  three example applications of the methods were presented. First, the wave scattering of the pentagonal object in free space was shown, demonstrating the use of the method for the problem with curvilinear boundary and free space boundary conditions. Second, a simulation with five round objects in free space was demonstrated. Finally the method was used to solve the inverse problem of locating a cylindrical underground body.

A future extension could be to improve the efficiency of the DG method used on the curvilinear body fitted grids by the use of an implicit timestepping method. This would allow the timesteps to be commensurate to those of the Hermite method at a relatively low cost since the linear systems needed to be inverted would be essentially one dimensional.  Another natural extension of this work would be to apply the techniques presented here to the  elastic wave equation.

\begin{acknowledgements}

This work was supported in part by the National Science Foundation under Grant NSF-1913076. Any opinions, findings, and conclusions or recommendations expressed in this material are those of the author(s) and do not necessarily reflect the views of the National Science Foundation.
\end{acknowledgements}
\section*{Conflict of interest}
The authors declare that they have no conflict of interest.
\bibliographystyle{spmpsci}
\bibliography{manuscript}
\end{document}